\documentclass[10pt,a4paper]{article}
\newif \ifSIAM

\newif \ifBW

\newif \ifELSEVIER

\usepackage{graphicx} 
\usepackage{epstopdf} 
\usepackage{amssymb}
\usepackage{amsmath} 
\usepackage{booktabs} 
\usepackage{paralist} 
\usepackage{bm} 
\usepackage{bbm} 
\usepackage{url} 
\usepackage{xcolor} 
\usepackage{mathtools} 
\usepackage{enumitem} 
\usepackage{placeins} 
\usepackage{float} 
\usepackage{subfig} 
\usepackage{algorithmicx} 
\usepackage{algorithm}
\usepackage{algpseudocode}
\usepackage{multirow}
\usepackage{siunitx}

\usepackage{ifpdf}

\algblockdefx[Module]{Module}{EndModule}%
[2][]{\textbf{module} #1(#2)}%
{\textbf{end module}}

\usepackage{amsthm}

\ifELSEVIER 
\numberwithin{equation}{section}
\else
\makeatletter
\providecommand\phantomcaption{\caption@refstepcounter\@captype}
\makeatother

\numberwithin{equation}{section}
\fi

\newif \iffract
\fracttrue

\newif \ifbbm

\newif \ifPDF
\PDFtrue

\newif \ifPdfTex
\PdfTextrue

\renewcommand{\algorithmicrequire}{$\rhd$ \textbf{Input:}}
\renewcommand{\algorithmicensure}{$\rhd$ \textbf{Output:}}

\definecolor{salem_green}{rgb}{0.1176,0.5098,0.2980}
\newcommand{\pd}[1]{{\color{black}{#1}}}

\definecolor{yellow}{rgb}{1, 0.8, 0}

\definecolor{ecstasy}{rgb}{0.9765,0.4118,0.0549}

\newif \ifwhere

\newcommand{\be}{\begin{equation}}
\newcommand{\ee}{\end{equation}}
\newcommand{\bea}{\begin{eqnarray}}
\newcommand{\eea}{\end{eqnarray}}
\newcommand{\bean}{\begin{eqnarray*}}
\newcommand{\eean}{\end{eqnarray*}}
\def\ba#1\ea{\begin{align}#1\end{align}}
\def\ban#1\ean{\begin{align*}#1\end{align*}}
\def\bat#1\eat{\begin{alignat}#1\end{alignat}}
\def\batn#1\eatn{\begin{alignat*}#1\end{alignat*}}
\def\bs#1\es{\begin{split}#1\end{split}}
\newcommand{\bse}{\begin{subequations}}
\newcommand{\ese}{\end{subequations}}
\newcommand{\bt}{\begin{theorem}}
\newcommand{\et}{\end{theorem}}
\newcommand{\bl}{\begin{lemma}}
\newcommand{\el}{\end{lemma}}
\newcommand{\bc}{\begin{corollary}}
\newcommand{\ec}{\end{corollary}}
\newcommand{\bp}{\begin{proof}}
\newcommand{\ep}{\end{proof}}
\newcommand{\bd}{\begin{definition}}
\newcommand{\ed}{\end{definition}}
\newcommand{\br}{\begin{remark}}
\newcommand{\er}{\end{remark}}
\newcommand{\bas}{\begin{assumption}}
\newcommand{\eas}{\end{assumption}}
\newcommand{\bex}{\begin{example}}
\newcommand{\eex}{\end{example}}
\newcommand{\bqo}{\begin{quote}}
\newcommand{\eqo}{\end{quote}}
\newcommand{\bdc}{\begin{description}}
\newcommand{\edc}{\end{description}}
\newcommand{\bi}{\begin{itemize}}
\newcommand{\ei}{\end{itemize}}
\newcommand{\ben}{\begin{enumerate}}
\newcommand{\een}{\end{enumerate}}
\newcommand{\bmar}{\begin{marking}}
\newcommand{\emar}{\end{marking}}
\newcommand{\bpro}{\begin{procedure}}
\newcommand{\epro}{\end{procedure}}

\ifELSEVIER
	\newtheorem{lemma}{Lemma}[section]
	\newtheorem{corollary}[lemma]{Corollary}
	\newtheorem{assumption}[lemma]{Assumption}
	\newtheorem{theorem}[lemma]{Theorem}
	\newtheorem{definition}[lemma]{Definition}
	\newtheorem{remark}[lemma]{Remark}
	\newtheorem{example}[lemma]{Example}
	\newtheorem{marking}[lemma]{Marking}
	\newtheorem{procedure}[lemma]{Procedure}
\else
\newtheorem{lemma}{Lemma}[section]
\newtheorem{corollary}[lemma]{Corollary}
\newtheorem{assumption}[lemma]{Assumption}
\newtheorem{theorem}[lemma]{Theorem}
\newtheorem{definition}[lemma]{Definition}
\newtheorem{remark}[lemma]{Remark}
\newtheorem{example}[lemma]{Example}
\newtheorem{marking}[lemma]{Marking}
\newtheorem{procedure}[lemma]{Procedure}
\fi

\newcommand\Om{\Omega}

\newcommand\omal{{\omega^\ver_\ell}}
\newcommand\omel{{\omega_\ell}}

\newcommand\Gr{\nabla}

\newcommand\Dv{\nabla {\cdot}}

\newcommand\dv{\mathrm{div}}

\newcommand\scp{{\cdot}}

\DeclarePairedDelimiter\norm{\lVert}{\rVert}

\ifbbm

\else

\fi

\newcommand\Hoo{H^1_0(\Om)}

\newcommand\Lt{L^2(\Om)}

\newcommand\Hdv{{\mathbf H}(\dv,\Om)}
\newcommand\Hdvi[1]{{\mathbf H}(\dv,#1)}

\ifbbm

\else

\fi

\newcommand\eg{e.g.}

\newcommand\eal{{\em et al.}}
\newcommand\eq{:=}

\newcommand\pt{\partial}

\newcommand\ra{\rightarrow}

\ifbbm

\else

\fi

\newcommand{\elm}{{K}}

\newcommand{\ver}{{\ta}}

\newcommand\Th{\mathcal{T}_\ell}
\newcommand\Thh{\mathcal{T}_{\ell+1}}

\newcommand\Ta{\mathcal{T}^{\ver}_\ell}
\newcommand\Tah{\mathcal{T}^{\ver,h}_\ell}
\newcommand\Tap{\mathcal{T}^{\ver,p}_\ell}
\newcommand\pa{\mathbf{p}^{\ver}_\ell}
\newcommand\pah{\mathbf{p}^{\ver,h}_\ell}
\newcommand\pap{\mathbf{p}^{\ver,p}_\ell}

\newcommand\Vh{\mathcal{V}_\ell}
\newcommand\Vhint{\mathcal{V}^{\mathrm{int}}_\ell}
\newcommand\Vhext{\mathcal{V}^{\mathrm{ext}}_\ell}
\newcommand\Vint[1]{\mathcal{V}^{\mathrm{int}}_{#1}}
\newcommand\Vext[1]{\mathcal{V}^{\mathrm{ext}}_{#1}}

\newcommand\VK{\mathcal{V}_\elm}

\newcommand\Vmarked{\pd{\widetilde{\V}_{\ell}^{\theta}}}
\newcommand\Vmarkedh{\widetilde{\V}_{\ell}^{h}}
\newcommand\Vmarkedp{\widetilde{\V}_{\ell}^{p}}

\newcommand\Mell{\pd{\M_{\ell}^{\theta}}}
\newcommand\Mellth{\M_{\ell}^{\theta_{\ell}}}
\newcommand\Mellh{{\M_{\ell}^{h}}}
\newcommand\Mellp{{\M_{\ell}^{p}}}

\newcommand\uu{u}
\newcommand\uh{u_\ell}
\newcommand{\uhp}{u_\ell^{p-1}}
\newcommand\uhh{u_{\ell+1}}

\newcommand\rahp{r^{\ver,hp}}
\newcommand\Vahp{V_\ell^{\ver,hp}}

\newcommand\frh{{\bm \sigma}_\ell}

\ifbbm

\else

\fi

\newcommand\Cred{C_\mathrm{red}}

\newcommand\ta{{\bf a}}

\newcommand\tn{{\bf n}}
\newcommand\tp{{\bf p}}

\newcommand\tv{{\bf v}}

\newcommand\tx{{\bf x}}

\newcommand\tV{{\bf V}}

\newcommand\M{\mathcal{M}}

\newcommand\T{\mathcal{T}}

\newcommand\V{\mathcal{V}}

\newcommand\RR{{\mathbb R}}

\ifbbm

\newcommand\PP{{\mathcal P}}

\else

\newcommand\PP{{\mathbb P}}

\fi

\iffract

\else

\fi

\ifBW 
\newcommand{\BW}{_BW}
\else
\newcommand{\BW}{}
\fi

\ifELSEVIER 
\else
	\title{An adaptive $hp$-refinement strategy with computable guaranteed bound on the error reduction factor\thanks{This
	project has received funding from the European Research Council (ERC)
	under the European Union's Horizon 2020 research and innovation program (grant agreement
	No 647134 GATIPOR).}}

	\author{Patrik Daniel\footnotemark[2] \footnotemark[3] \and Alexandre Ern\footnotemark[3] \footnotemark[2] \and Iain Smears\footnotemark[4] \and Martin Vohral\'ik\footnotemark[2] \footnotemark[3]}
\fi

\ifPdfTex \usepackage[breaklinks,bookmarks=false]{hyperref}
\else
\ifPDF \usepackage[dvips,breaklinks,bookmarks=false]{hyperref} \fi \fi

\ifPDF
\hypersetup{colorlinks, linkcolor=blue, citecolor=blue,
urlcolor=blue, plainpages=false, pdfwindowui=false,
pdfstartview={FitH}, pdftitle={An adaptive hp-refinement strategy with computable guaranteed bound on the error reduction factor}, pdfauthor={Patrik Daniel}}
\fi

\begin{document}

\ifELSEVIER \begin{frontmatter}
\fi

\ifELSEVIER
\else
	\maketitle

	\renewcommand{\thefootnote}{\fnsymbol{footnote}}

	\footnotetext[2]{Inria, 2 rue Simone Iff, 75589 Paris, France
	}

	\footnotetext[3]{
	Université Paris-Est, CERMICS (ENPC), 77455 Marne-la-Vallée 2, France
	}

	\footnotetext[4]{Department of Mathematics, University College London, Gower Street, WC1E 6BT London, United Kingdom
	}

	\renewcommand{\thefootnote}{\arabic{footnote}}
\fi

\begin{abstract}
We propose a new practical adaptive refinement strategy for $hp$-finite
element approximations of elliptic problems. Following recent theoretical
developments in polynomial-degree-robust a posteriori error analysis, we
solve two types of discrete local problems on vertex-based patches. The
first type involves the solution on each patch of a mixed finite element
problem with homogeneous Neumann boundary conditions, which leads to an
$\Hdv$-conforming equilibrated flux. This, in turn, yields a guaranteed
upper bound on the error and serves to mark mesh vertices for refinement
via D\"orfler's bulk-chasing criterion. The second type of local problems
involves the solution, on patches associated with
marked vertices only, of two separate primal finite element problems with
homogeneous Dirichlet boundary conditions, which serve to decide between 
\mbox{$h$-,} \mbox{$p$-,} or $hp$-refinement. Altogether, we show that these ingredients
lead to a computable guaranteed bound on the ratio of the errors between
successive refinements (error reduction factor). In a series of
numerical experiments featuring smooth and singular solutions, we study the
performance of the proposed $hp$-adaptive strategy and observe exponential
convergence rates. We also investigate the accuracy of
our bound on the reduction factor by evaluating the
ratio of the predicted reduction factor relative to the true error
reduction, and we find that this ratio is in general quite close to the
optimal value of one.
\end{abstract}

\ifELSEVIER
\else
	\bigskip

\fi


\ifELSEVIER
	\end{frontmatter}
	\linenumbers
\else
\fi

\section{Introduction} 

Adaptive discretization methods constitute an important tool in computational
science and engineering. Since the pioneering works on the $hp$-finite
element method by Gui and Babu\v{s}ka~\cite{Gui_Bab_hp_II_86,
Gui_Bab_hp_III_86} and Babu\v{s}ka and
Guo~\cite{Bab_Guo_hp_elem_1_86,Bab_Guo_hp_elem_2_86} in the 1980s, where it
was shown that for one-dimensional problems $hp$-refinement leads to
exponential convergence with respect to the number of degrees of freedom on
{\em a priori} adapted meshes, there has been a great amount of work devoted
to developing {\em adaptive $hp$-refinement} strategies based on {\em a
posteriori error estimates}. Convergence of $hp$-adaptive finite element
approximations for elliptic problems, has, though, been addressed only very
recently in~D{\"o}rfler and Heuveline~\cite{Dorf_Heu_hp_cvg_1D_07}, B{\"u}rg
and D{\"o}rfler~\cite{Burg_Dorf_hp_cvg_11}, and Bank, Parsania, and
Sauter~\cite{Bank_Pars_Saut_satur_est_hpFEM_13}. The first optimality result
we are aware of is by Canuto~\eal\ \cite{Can_Noch_Stev_Ver_hp_AFEM_16}, where
an important ingredient is the $hp$-coarsening routine by
Binev~\cite{Bin_hp_approx_13,Bin_hp_adap_15}. These works extend to the
$hp$-context the previous $h$-convergence and optimality results by
D\"orfler~\cite{Dorfler_marking_96}, Morin, Nochetto, and Siebert
\cite{Morin_converge_afem_02, Morin_loc_prob_stars_03},
Stevenson~\cite{Stev_opt_FE_07}, Casc{\'o}n~\eal\
\cite{Casc_Kreu_Noch_Sieb_08}, Carstensen~\eal\
\cite{Cars_Feis_Page_Prea_ax_adpt_14}, see also Nochetto {\em et
al.}~\cite{Noch_Sieb_Vees_09} and the references therein. It is worth
mentioning that most of the available convergence results are formulated for
adaptive methods driven by residual-type a posteriori error estimators; other
estimators have in particular been addressed in Casc\'on and
Nochetto~\cite{Casc_Noch_cvg_non_res_11} and Kreuzer and
Siebert~\cite{Kreu_Sieb_decay_afem_11}.

A key ingredient for adaptive $hp$-refinement is a local criterion in each
mesh cell marked for refinement that allows one to decide whether $h$-, $p$-,
or $hp$-refinement should be performed. There is a substantial amount of such
criteria proposed in the literature; a computational overview can be found in
Mitchell and McClain~\cite{Mitchell_hp_review_14,Mitchell_hp_review_long}. Some of the mathematically
motivated $hp$-decision criteria include, among others, those
proposed by Eibner and Melenk~\cite{Eib_Mel_hp_anal_07}, Houston and
S\"uli~\cite{Hous_Sul_hp_adpt_05} which both estimate the local regularity of the
exact weak solution. Our proposed strategy fits into the group of algorithms
based on solving local boundary value problems allowing us to forecast the
benefits of performing $h$- or $p$-refinement, as recently considered in,
\eg,~\cite{Burg_Dorf_hp_cvg_11,Dorf_Heu_hp_cvg_1D_07}. Similarly
to~\cite{Dorf_Heu_hp_cvg_1D_07}, we use the \emph{local} finite element
spaces associated with a specific type of refinement to perform the above
forecast and to take the local $hp$-refinement decision. We also mention the
work of Demkowicz \eal\ \cite{Dem_full_aut_hp_02} for an earlier, yet more
expensive, version of the look-ahead idea, where it is proposed to solve an
auxiliary problem on a \emph{global} finite element space corresponding to a
mesh refined uniformly either in $h$ or in $p$.

In the present work, we focus on the Poisson model problem with (homogeneous) Dirichlet
boundary conditions.
In weak form, the model problem reads as
follows: Find $\uu \in \Hoo$ such that
\be
	\label{eq_weak_form}
	\left( \Gr \uu , \Gr v \right) = (f,v) \qquad \forall v \in \Hoo,
\ee
where $\Om \subset \RR^d$, $d=2,3$, is a polygonal/polyhedral domain (open,
bounded, and connected set) with a Lipschitz boundary, $f \in \Lt$, $\Hoo$
denotes the Sobolev space of all functions in $\Lt$ which have all their
first-order weak derivatives in $\Lt$ and a zero trace on $\pt \Om$, and
$(\cdot, \cdot)$ stands for the $\Lt$ or $[\Lt]^d$~inner product. Our first
goal is to propose a reliable and computationally-efficient
\emph{$hp$-adaptive strategy} to approximate the model
problem~\eqref{eq_weak_form} that hinges on the recent theoretical
developments on polynomial-degree-robust a posteriori error estimates due to
Braess~\eal\ \cite{Brae_Pill_Sch_p_rob_09} and Ern and
Vohral\'ik~\cite{Ern_Voh_p_rob_15, Ern_Voh_p_rob_3D_16}. The present
$hp$-adaptive algorithm follows the well-established paradigm based
on an iterative loop where each step consists of the following four
modules:
\be	
	\label{eq_afem_loop}
	\mathrm{SOLVE} \ra \mathrm{ESTIMATE} \ra \mathrm{MARK} \ra 	
    \mathrm{REFINE}.
\ee
Here, $\mathrm{SOLVE}$ stands for application of the conforming finite
element method on a matching (no hanging nodes) simplicial mesh to
approximate the model problem~\eqref{eq_weak_form}; spatially-varying
polynomial degree is allowed. The module $\mathrm{ESTIMATE}$ is based on an
equilibrated flux a posteriori error estimate, obtained by solving, for each
mesh vertex, a local mixed finite element problem with a (homogeneous) Neumann
boundary condition on the patch of cells sharing the given vertex. The module
$\mathrm{MARK}$ is based on a bulk-chasing criterion inspired by the
well-known D\"orfler's marking~\cite{Dorfler_marking_96}; here we mark mesh
vertices and not simplices since we observe a smoother performance in
practice and since we later work with some vertex-based auxiliary
quantities.

The module $\mathrm{REFINE}$, where we include our $hp$-decision
criterion, is organized into three steps. First, we solve two local finite
element problems on each patch of simplices attached to a mesh vertex marked
for refinement, with either the mesh refined or the polynomial degree
increased. This is inspired by the key observation
from~\cite[Lemma~3.23]{Ern_Voh_p_rob_15} that guaranteed local efficiency can
be materialized by some local conforming finite element solves. These
conforming residual liftings allow us, in particular, to estimate the effect
of applying $h$- or $p$-refinement, and lead to a partition of the set of
marked vertices into two disjoint subsets, one collecting the mesh vertices
flagged for $h$-refinement and the other collecting the mesh vertices
flagged for $p$-refinement. The second step of the module
$\mathrm{REFINE}$ uses these two subsets to flag the simplices for $h$-,
$p$, or $hp$-refinement.
Finally, the third step of the module $\mathrm{REFINE}$ uses
the above sets of flagged simplices to build the next simplicial mesh
and the next polynomial-degree distribution. Let us mention that recently,
Dolej{\v{s}}{\'{\i}} \eal\ \cite{hp_Dol_Ern_Voh_16} also devised an
$hp$-adaptive algorithm driven by polynomial-degree-robust a posteriori error
estimates based on the equilibrated fluxes from
\cite{Brae_Pill_Sch_p_rob_09,Ern_Voh_p_rob_15,Ern_Voh_p_rob_3D_16}. The
differences with the present work are that the interior penalty discontinuous
Galerkin method is considered in \cite{hp_Dol_Ern_Voh_16}, and more
importantly, that the present $hp$-decision criterion hinges on local primal
solves on patches around marked vertices.

The second goal of the present work is to show that the proposed
$hp$-adaptive strategy automatically leads to a {\em computable guaranteed}
{\em bound} on the {\em error reduction factor} between two
consecutive steps of the adaptive loop~\eqref{eq_afem_loop}. More precisely,
we show how to compute explicitly a real number $\Cred\in [0,1]$ so that
\begin{equation} \label{eq_C_red}
  	\norm{\Gr (\uu -\uhh)} \le \Cred \norm{\Gr (\uu - \uh)},
\end{equation}
where $\uh$ denotes the discrete solution on $\ell$-th iteration of the
adaptive loop, see Theorem~\ref{th_contraction} below. Thus the number
$\Cred$ gives a guaranteed (constant-free) bound on the ratio of the errors
between successive refinements. This must not be confused with saying that
the error is guaranteed to be reduced, since the case $\Cred=1$ cannot be
ruled out in general without additional assumptions (e.g.\ an interior node
property, see \cite{Morin_converge_afem_02} for further details). The
computation of $\Cred$ crucially relies on a combined use of the equilibrated
fluxes and of the conforming residual liftings, which were already
 used for the error estimation and $hp$-refinement decision
criterion respectively. It is worth noting that we consider a homogeneous Dirichlet
boundary condition for the local residual liftings in order to obtain an
estimate on the error reduction factor that is as sharp as possible.

The rest of this manuscript is organized as follows.
Section~\ref{sec:setting} describes the discrete setting and introduces some
useful notation. Section~\ref{sec_SolEstMark} presents the modules
$\mathrm{SOLVE}$, $\mathrm{ESTIMATE}$, and $\mathrm{MARK}$, whereas
Section~\ref{sec_hp_strat} presents the module $\mathrm{REFINE}$.
Section~\ref{sec_err_red} contains our main result on a guaranteed bound on
the error reduction factor. Finally, numerical experiments on
two-dimensional test cases featuring smooth and singular solutions are
discussed in Section~\ref{sec_num}, and conclusions are drawn in
Section~\ref{sec:conc}.

\section{Discrete setting}
\label{sec:setting}

The main purpose of the adaptive loop~\eqref{eq_afem_loop} is to generate a
sequence of finite-dimensional $H^1_0$-conforming finite element spaces
$(V_\ell)_{\ell\geq0}$, where the integer $\ell\geq0$ stands for the
iteration counter in~\eqref{eq_afem_loop}. $H^1_0$-conformity means that
$V_\ell\subset H^1_0(\Omega)$ for all $\ell\geq0$. In this work, we shall
make the following nestedness assumption:
\begin{equation} \label{eq:nested}
V_\ell \subset V_{\ell+1}, \qquad \forall \ell \geq 0.
\end{equation}
The space $V_\ell$ is built from two ingredients: \textup{(i)}
a matching simplicial mesh $\Th$ of the computational domain $\Om$, that is, a finite collection of (closed) non-overlapping simplices $\elm\in\Th$ covering $\overline\Om$ exactly and such that the intersection of two different simplices is
either empty, a common vertex, a common edge, or a common face;
\textup{(ii)} a polynomial-degree distribution described by the
vector $\tp_\ell \eq (p_{\ell,\elm})_{\elm \in \Th}$ that assigns
a polynomial degree to each simplex $\elm \in \Th$.
The conforming finite element space $V_\ell$ is then defined as
\begin{equation}
	\nonumber
	V_\ell \eq \PP_{\tp_\ell}(\Th) \cap \Hoo, \qquad \forall \ell \geq 0,
\end{equation}
where $\PP_{\tp_\ell}(\Th)$ denotes the space of piecewise polynomials of
total  degree $p_{\ell,\elm} \geq 1$ on each simplex $\elm \in \Th$. In
other words, any function $v_\ell\in V_\ell$ satisfies $v_\ell\in
H^1_0(\Omega)$ and $v_{\ell}|_\elm \in \PP_{p_{\ell,K}}(\elm)$ for all $\elm
\in \Th$, where for an integer $p\geq 1$, $\PP_p (\elm)$ stands for the space
of polynomials of total degree at most $p$ on the simplex $\elm$.

The initial mesh $\T_0$ and the initial polynomial-degree distribution $\tp_0$ are given, and the purpose of each step $\ell\ge0$ of the adaptive
loop~\eqref{eq_afem_loop} is to produce the next mesh $\Thh$ and the next polynomial-degree distribution $\tp_{\ell+1}$. In order to ensure the nestedness property~\eqref{eq:nested}, the following two properties are to be satisfied: \textup{(i)} The sequence $(\Th)_{\ell \geq 0}$ is hierarchical, i.e., for all $\ell\geq0$ and all $\widetilde{\elm}\in \Thh$, there is a unique
simplex $\elm \in \Th$, called the parent of $\widetilde{\elm}$
so that $\widetilde{\elm}\subseteq K$; \textup{(ii)}
The local polynomial degree is locally increasing, i.e.,
for all $\ell\geq0$ and all $\widetilde{\elm}\in \Thh$,
$p_{\ell+1,\widetilde{\elm}}\geq p_{\ell,\elm}$ where $\elm \in \Th$
is the parent of $\widetilde{\elm}$.
Moreover, we assume the following shape-regularity property: There exists a
constant $\kappa_{\T}>0$ such that $\max_{\elm\in\Th} h_\elm / \rho_\elm \le
\kappa_{\T}$ for all $\ell \geq 0$, where $h_\elm$ is the diameter of
$\elm$ and $\rho_\elm$ is the diameter of the largest ball inscribed
in $\elm$.

Before closing this section, we introduce some further useful notation. The
set of vertices $\Vh$ of each mesh $\Th$ is decomposed into $\Vhint$ and
$\Vhext$, the set of inner and boundary vertices, respectively.
For each vertex $\ver \in \Vh$, the so-called hat function~$\psi^\ver_\ell$ is the continuous, piecewise affine function that takes the value $1$ at the vertex $\ver$ and the value $0$ at all the other vertices of $\V_\ell$; the function $\psi^\ver_\ell$ is in $V_\ell$ for all $\ver\in\Vhint$.
Moreover, we consider the simplex patch $\Ta\subset \Th$ which is the collection of the simplices in $\Th$ sharing the vertex
$\ver \in \Vh$, and we denote by $\omal$ the corresponding open subdomain. Finally, for
each simplex $\elm \in \Th$, $\V_\elm$ denotes the set of vertices of $\elm$.

\section{The modules SOLVE, ESTIMATE, and MARK} \label{sec_SolEstMark}

In this section we present the modules SOLVE, ESTIMATE, and MARK from
the adaptive loop~\eqref{eq_afem_loop}. Let $\ell\geq0$ denote the current
iteration number.

\subsection{The module SOLVE} 

The module SOLVE takes as input the $H^1_0$-conforming finite element space $V_\ell$ and outputs the discrete function $\uh \in V_\ell$ which is the unique solution of
\be \label{eq_fem_model}
	(\Gr \uh, \Gr v_\ell) = (f, v_\ell) \qquad \forall v_\ell \in V_\ell.
\ee

\subsection{The module ESTIMATE}

Following~\cite{Dest_Met_expl_err_CFE_99,Brae_Pill_Sch_p_rob_09,Ern_Voh_p_rob_15,
hp_Dol_Ern_Voh_16, Ern_Voh_p_rob_3D_16}, see also the references therein, the
module ESTIMATE relies on an equilibrated flux a posteriori error estimate on
the energy error $\norm{\Gr(\uu - \uh)}$. The module ESTIMATE takes as input
the finite element solution $u_\ell$ and outputs a collection of local error
indicators $\{\eta_\elm\}_{\elm\in \Th}$. The equilibrated flux is constructed locally from
mixed finite element solves on the simplex patches $\Ta$ attached to each
vertex $\ver\in\Vh$. For this construction, we consider as
in~\cite{hp_Dol_Ern_Voh_16} the local polynomial degree
$p_\ver^{\mathrm{est}} \eq \max_{\elm \in \Ta} p_{\ell,\elm}$ (any other
choice so that $p_\ver^{\mathrm{est}} \geq \max_{\elm \in \Ta} p_{\ell,\elm}$
can also be employed). We consider the local
Raviart--Thomas--N\'ed\'elec mixed finite element spaces
$(\tV_\ell^\ver,Q_\ell^\ver)$ 
\ifELSEVIER
\else
which are defined for all $\ver \in \Vint{\ell}$ by
    \ban
    \tV_\ell^\ver  &\eq  \{ \tv_\ell \in \Hdvi{\omal}; \,
    \tv_\ell|_\elm \in \mathbf{RTN}_{p_\ver^{\mathrm{est}}}(\elm), \,
    \forall \elm \in \Ta, \, \tv_\ell \scp \tn_{\omal} = 0 \text{ on } \pt \omal\}, \\[1mm]
       Q_\ell^\ver  &\eq  \{ q_\ell \in \PP_{p_\ver^{\mathrm{est}}}(\Ta); \, (q_\ell,1)_{\omal} = 0 \},
    \ean
and, for all $\ver \in \Vext{\ell}$,
    \ban    
    \tV_\ell^\ver  \!&\eq\! \{ \tv_\ell \! \in\! \Hdvi{\omal}; \,
    \tv_\ell|_\elm\! \in\!  \mathbf{RTN}_{p_\ver^{\mathrm{est}}}(\elm),\,
    \forall \!\elm\! \in\! \Ta,\, \tv_\ell \scp \tn_{\omal} \!= 0 \text{ on } \pt \omal\! \setminus\! \pt \Om\}, \\
    Q_\ell^\ver \!&\eq\! \PP_{p_\ver^{\mathrm{est}}}(\Ta),    
    \ean
where $\mathbf{RTN}_{p_\ver^{\mathrm{est}}}(\elm) \eq [\PP_{p_\ver^{\mathrm{est}}}(\elm)]^d + \PP_{p_\ver^{\mathrm{est}}}(\elm) \tx$, and $\tn_{\omal}$ denotes the unit outward-pointing normal to $\omal$.
\fi

\bd[Flux reconstruction $\frh$]\label{def_equi_flux_construction}
Let $u_\ell$ solve~\eqref{eq_fem_model}. The
global equilibrated flux $\frh$ is constructed as $\frh
\eq \sum_{\ver \in \Vh} \frh^\ver$, where, for each vertex $\ver \in \Vh$,
$(\frh^\ver, \gamma^\ver_\ell) \in \tV_\ell^\ver \times Q_\ell^\ver$ solves
\bse \label{eq_dual_mixed_problem}
\bat{2} 
	(\frh^\ver, \tv_\ell)_{\omal} - (\gamma^\ver_\ell, \Dv \tv_\ell)_{\omal} &=-(\psi^\ver_\ell \Gr \uh, \tv_\ell)_{\omal} \quad &\forall \tv_\ell \in \tV_\ell^\ver ,\nonumber\\
	(\Dv \frh^\ver, q_\ell)_{\omal} &=(f \psi^\ver_\ell  - \Gr \uh \scp \Gr
\psi^\ver_\ell,q_\ell)_{\omal} \quad &\forall q_\ell \in Q_\ell^\ver;\nonumber \eat \ese
or, equivalently,
\be
	\nonumber
    \frh^\ver \eq \arg\min_{\tv_\ell \in \tV_\ell^\ver, \, \Dv \tv_\ell =
    \Pi_{Q_\ell^\ver} (f \psi^\ver_\ell - \Gr \uh \scp \Gr
    \psi^\ver_\ell)} \norm{\psi^\ver_\ell \Gr \uh + \tv_\ell}_{\omal},
\ee
and where $\frh^\ver$ is extended by zero outside $\omal$.
\ed

Note that the Neumann compatibility condition for the problem
\eqref{eq_dual_mixed_problem} is satisfied for all $\ver\in\Vhint$ (take $v_\ell = \psi^\ver_\ell$ as a test function in~\eqref{eq_fem_model}). Moreover,
Definition~\ref{def_equi_flux_construction} yields a globally
$\Hdv$-conforming flux reconstruction $\frh$ such that, for all $\elm \in
\Th$, $(\Dv \frh, v_\ell)_\elm = (f, v_\ell)_\elm$ for all $v_\ell \in
\PP_{\min_{\ver \in \VK} p_\ver^{\mathrm{est}}}(\elm)$,
see~\cite[Lemma~3.6]{hp_Dol_Ern_Voh_16}. Using the current
notation,~\cite[Theorem~3.3]{hp_Dol_Ern_Voh_16} states the following result.
\bt[Guaranteed upper bound on the error] \label{thm_est} Let $\uu$ solve \eqref{eq_weak_form} and $\uh$ solve \eqref{eq_fem_model}. Let $\frh$ be the equilibrated
flux reconstruction of Definition~\ref{def_equi_flux_construction}. Then
\be \label{eq_est}
    \norm{\Gr(\uu - \uh)} \leq \eta({\Th}) \eq \left\lbrace \sum_{\elm \in \Th} \eta_\elm^2
    \right\rbrace^\frac{1}{2}, \qquad \eta_\elm \eq \norm{\Gr \uh + \frh}_\elm + \frac{h_\elm}{\pi}\norm{f - \Dv \frh}_\elm.
\ee
\et

As discussed in, \eg, \cite[Remark~3.6]{Ern_Voh_p_rob_15}, the term
$\frac{h_\elm}{\pi}\norm{f - \Dv \frh}_\elm$ represents, for all $\elm\in\Th$,
a local oscillation
in the source datum $f$ that, under suitable smoothness assumptions,
converges to zero two orders faster than the
error. To cover the whole computational range in our numerical experiments, this term is kept in the error indicator $\eta_\elm$.

\subsection{The module MARK} 

The module MARK takes as input the local error estimators
$\{\eta_\elm\}_{\elm\in\Th}$ from Theorem~\ref{thm_est} and outputs a set of
marked vertices $\Vmarked\subset \Vh$ using a bulk-chasing criterion inspired
by the well-known D\"orfler's marking criterion~\cite{Dorfler_marking_96}.
The reason why we mark vertices and not simplices is that our $hp$-decision
criterion in the module REFINE (see Section~\ref{sec_hp_strat} below) hinges
on the solution of local primal solves posed on the patches $\Ta$ associated with
the marked vertices $\ver \in \Vmarked$; we also observe in practice a smoother
performance of the overall $hp$-adaptive procedure when marking vertices than when marking 
elements. Vertex-marking strategies are also considered, among others,
in~\cite{Morin_converge_afem_02,Can_Noch_Stev_Ver_sat_16}.

For a subset $\mathcal{S} \subset \Th$, we use the notation
$\eta(\mathcal{S}) \eq \lbrace \sum_{\elm \in \mathcal{S}} \eta_\elm^2
\rbrace^{1/2}$. In the module MARK, the set of marked vertices $\Vmarked$ is
selected in such a way that
\be \label{eq:def_Vmarked}
    \eta\bigg(\bigcup_{\ver \in \Vmarked} \Ta\bigg) \geq \theta\,\eta({\Th}),
\ee
where $\theta \in (0,1]$ is a fixed threshold. Letting
\be \label{eq:def_Ml}
    \Mell\eq \bigcup_{\ver \in \Vmarked} \Ta \subset \Th
\ee
be the collection of all the simplices that belong to a patch associated
with a marked vertex, we observe that \eqref{eq:def_Vmarked} means that
$\eta(\Mell) \geq \theta\,\eta({\Th})$. To select a set $\Vmarked$ of
minimal cardinality, the mesh vertices in $\Vh$ are sorted by comparing the
vertex-based error estimators $\eta(\Ta)$ for all $\ver\in \Vh$, and a greedy
algorithm is employed to build the set $\Vmarked$. The module MARK is
summarized in Algorithm~\ref{def_mark_patch_alg}. A possibly slightly larger
set $\Vmarked$ can be constructed with linear cost in terms of the number of
mesh vertices by using the algorithm proposed
in~\cite[Section~5.2]{Dorfler_marking_96}.

\begin{algorithm}
\caption{(\textbf{module MARK})}\label{def_mark_patch_alg}
\begin{algorithmic}[1]
\Procedure{MARK}{$\{\eta_\elm\}_{\elm \in \Th}$, $\theta$}
 \State \algorithmicrequire{ error indicators $\{\eta_\elm\}_{\elm \in \Th}$, marking parameter $\theta\in (0,1]$}
 \State \algorithmicensure{ set of marked vertices $\Vmarked$}
\ForAll{$\ver \in \V_\ell$}
	\State Compute the vertex-based error estimator $\eta(\Ta)$
\EndFor

\State Sort the vertices according to $\eta(\Ta)$
\State Set 
	   $\Vmarked \eq \emptyset$

\While{${\eta(\bigcup_{\ver \in \Vmarked} \Ta)} < \theta\,\eta({\Th})$}
\State {Add to $\Vmarked$ the next sorted vertex $\ver\in \Vh\setminus \Vmarked$}
\EndWhile
\EndProcedure
\end{algorithmic}
\end{algorithm}
\FloatBarrier

\section{The module REFINE}
\label{sec_hp_strat}

The module REFINE takes as input the set of marked vertices $\Vmarked$ and
outputs the mesh $\Thh$ and the polynomial-degree distribution $\tp_{\ell+1}$
to be used at the next step of the adaptive loop~\eqref{eq_afem_loop};
the integer $\ell\geq0$ is the current iteration number therein. This
module is organized into three steps. First, an $hp$-decision is made on all
the marked vertices so that each marked vertex $\ver\in\Vmarked$ is flagged
either for $h$-refinement or for $p$-refinement. This means that the set
$\Vmarked$ is split into two disjoint subsets $\Vmarked=\Vmarkedh \cup
\Vmarkedp$ with obvious notation (here we drop the superscript $\theta$ 
to simplify the notation). Then, in the second step, the subsets
$\Vmarkedh$ and $\Vmarkedp$ are used to define subsets $\Mellh$ and
$\Mellp$ of the set of marked simplices $\Mell$ (see~\eqref{eq:def_Ml}).
The subsets $\Mellh$ and $\Mellp$ are not necessarily disjoint which
means that some simplices can be flagged for $hp$-refinement. Finally, the
two subsets $\Mellh$ and $\Mellp$ are used to construct $\Thh$ and
$\tp_{\ell+1}$.

\subsection{$hp$-decision on vertices}

Our $hp$-decision on marked vertices is made on the basis of two local primal solves on the patch $\Ta$ attached to each
marked vertex $\ver\in \Vmarked$. The idea is to construct two distinct local
patch-based spaces in order to emulate separately the effects of $h$- and
$p$-refinement. Let us denote the polynomial-degree distribution in the patch
$\Ta$ by the vector $\pa \eq (p_{\ell,\elm})_{\elm\in\Ta}$.

\begin{figure}[!h]
\begin{minipage}[b]{0.295\linewidth}
	\centering
	\includegraphics[width=\textwidth]{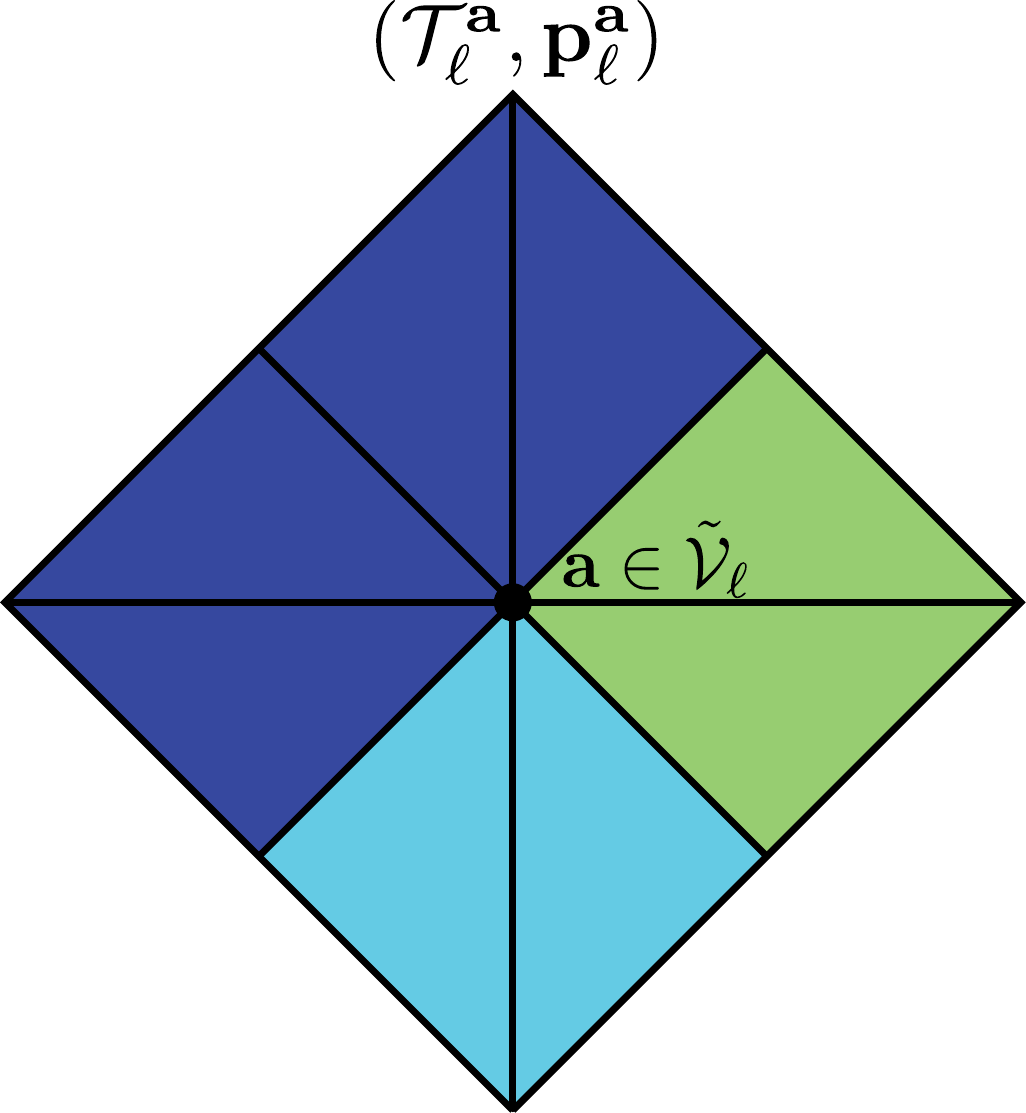}
\end{minipage}
\hfill
\begin{minipage}[b]{0.295\linewidth}
	\centering
	\includegraphics[width=\textwidth]{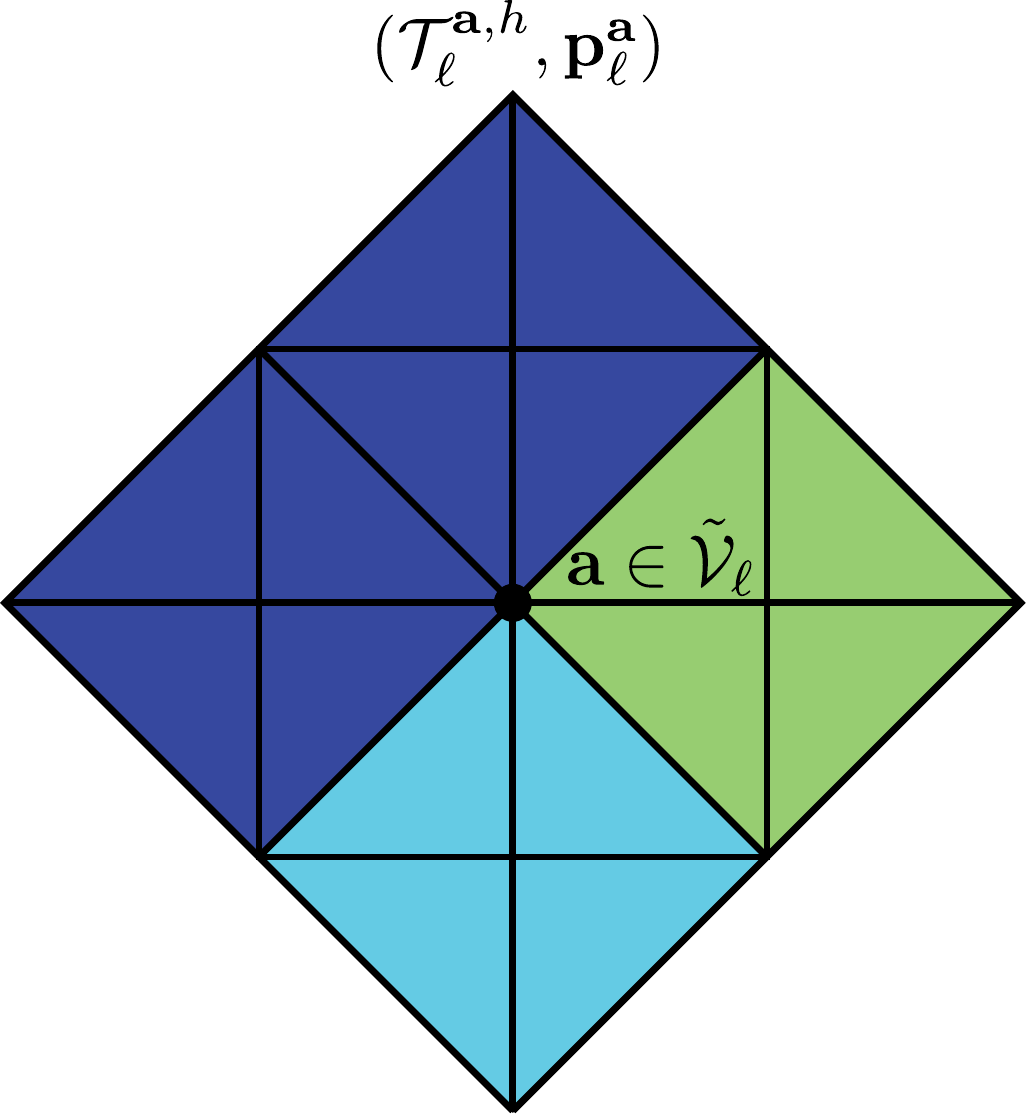}
\end{minipage}
\hfill
\begin{minipage}[b]{0.295\linewidth}
	\centering
	\includegraphics[width=\textwidth]{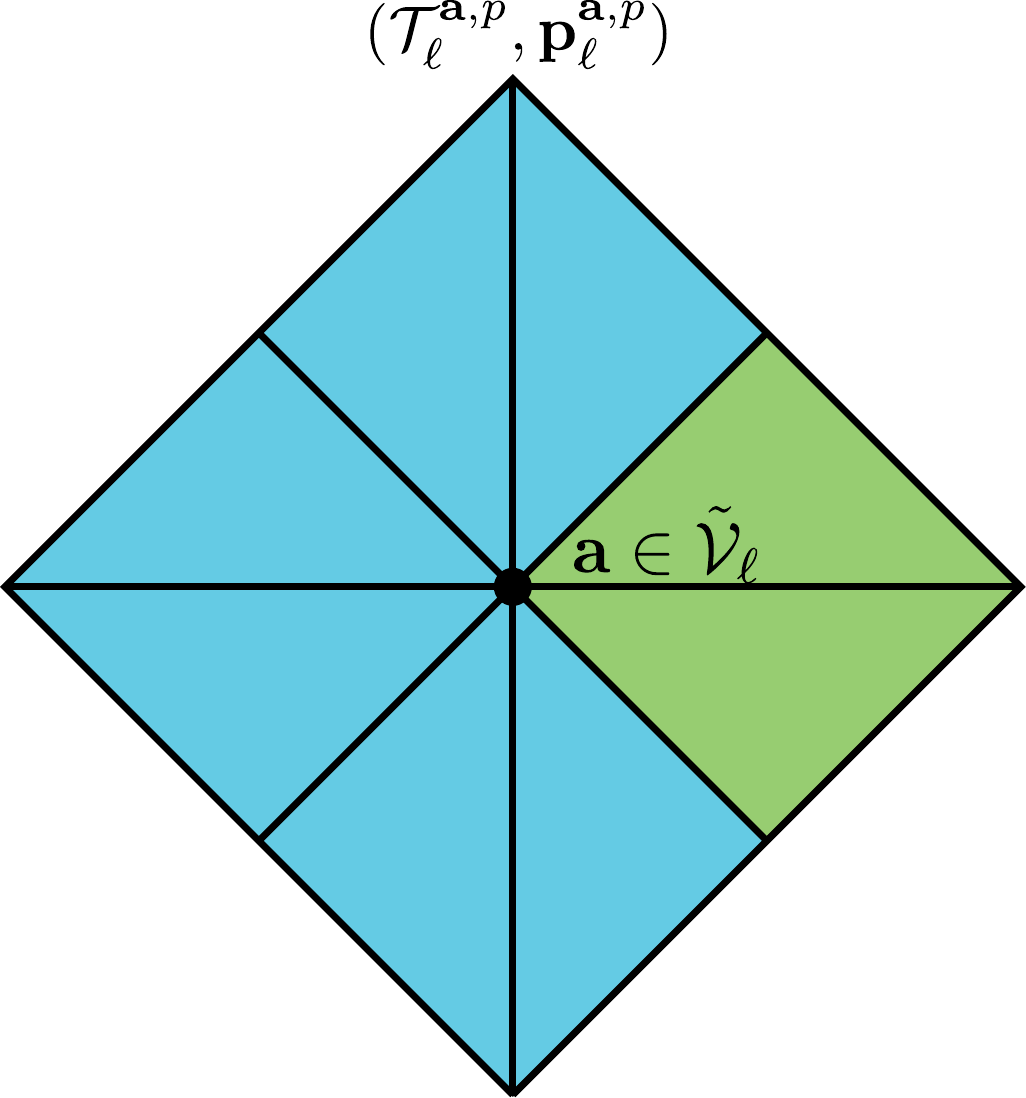}
\end{minipage}
\hfill
\begin{minipage}[b]{0.05\linewidth}
	\centering		
	\includegraphics[width=\textwidth]{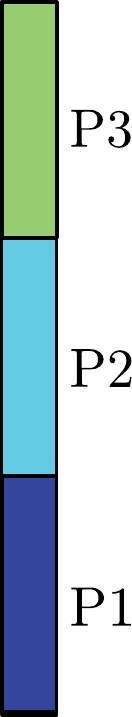}
\end{minipage}
	\caption{An example of patch $\T_\ell^\ver$ together with its polynomial-degree
    distribution $\tp_{\ell}^{\ver}$ (\textit{left}) and its $h$-refined (\textit{center})
    and $p$-refined versions (\textit{right}) from Definitions~\ref{def_h_ref} and \ref{def_p_ref} respectively.}
\end{figure}
\bd[\textbf{$h$-refinement residual}] \label{def_h_ref} Let $\ver\in\Vmarked$
be a marked vertex with associated patch $\Ta$ and polynomial-degree
distribution $\tp_\ell^\ver$. We set
\be \label{eq_V_a_h} V_\ell^{\ver,h} \coloneqq \PP_{\pah}(\Tah)\cap
H^1_0(\omal), \ee
where $\Tah$ is obtained as a matching simplicial refinement of $\Ta$ by dividing each simplex $\elm \in \Ta$ into at least two children simplices, and the polynomial-degree distribution $\pah$ is obtained from $\pa$ by assigning to each newly-created simplex the same polynomial degree as its parent. Then, we let $r^{\ver,h} \in V_\ell^{\ver,h}$ solve
\be
\nonumber
(\Gr r^{\ver,h}, \Gr v^{\ver,h})_{\omal} = (f, v^{\ver,h})_{\omal}-(\Gr u_{\ell}, \Gr v^{\ver,h})_{\omal} \qquad\forall\, v^{\ver,h} \in V_\ell^{\ver,h}.
\ee
\ed

\bd[\textbf{$p$-refinement residual}] \label{def_p_ref} Let $\ver\in\Vmarked$
be a marked vertex with associated patch $\Ta$ and polynomial-degree
distribution $\tp_\ell^\ver$. We set
\be \label{eq_V_a_p} V^{\ver,p}_\ell \coloneqq \PP_{\pap}(\Tap)\cap
H^1_0(\omal), \ee
where $\Tap \eq \Ta$, and the polynomial-degree distribution $\pap$ is obtained from $\pa$ by assigning to each simplex $K\in \Tap=\Ta$ the polynomial degree $p_{\ell,K}+\delta^\ver_\elm$ where
\be \label{eq_increm_p} \delta^\ver_\elm \eq \left\{ \begin{array}{ll}
		1  & \mbox{\textup{if}}\; p_{\ell,\elm} = \min_{\elm' \in \Ta} p_{\ell,\elm'}, \\
		0 & \mbox{\textup{otherwise}}.
\end{array}
\right.
\ee
Then, we let $r^{\ver,p} \in V^{\ver,p}_\ell$ solve
\be
\nonumber
(\Gr r^{\ver,p}, \Gr v^{\ver,p})_{\omal} = (f, v^{\ver,p})_{\omal}-(\Gr u_{\ell}, \Gr v^{\ver,p})_{\omal}  \qquad \forall\, v^{\ver,p} \in V^{\ver,p}_\ell.
\ee
\ed

The local residual liftings $r^{\ver,h}$ and $r^{\ver,p}$ from
Definitions~\ref{def_h_ref} and~\ref{def_p_ref}, respectively, are used to
define the following two disjoint subsets of the set of marked vertices
$\Vmarked$:
\bse \label{eq:decision_vertex}
\begin{align}
\Vmarkedh &\eq \{ \ver\in\Vmarked\;|\; \norm{\Gr r^{\ver,h}}_{\omal} \geq \norm{\Gr r^{\ver,p}}_{\omal} \},\\
\Vmarkedp &\eq \{ \ver\in\Vmarked\;|\; \norm{\Gr r^{\ver,h}}_{\omal} < \norm{\Gr r^{\ver,p}}_{\omal} \},
\end{align}
\ese
in such a way that
\be
\nonumber
\Vmarked = \Vmarkedh \cup \Vmarkedp, \qquad \Vmarkedh \cap \Vmarkedp = \emptyset.
\ee
The above $hp$-decision criterion on vertices means that a marked vertex is
flagged for $h$-refinement if the local residual norm $\norm{\Gr
r^{\ver,h}}_{\omal}$ is larger than $\norm{\Gr r^{\ver,p}}_{\omal}$;
otherwise, this vertex is flagged for $p$-refinement. Further motivation
for this choice is discussed in Remark~\ref{rem_loc_res} below.

\br[$p$-refinement] Other choices are possible for the polynomial-degree
increment defined in~\eqref{eq_increm_p}. One possibility is to set
$\delta^\ver_\elm=1$ for all $\elm\in \Ta$. However, in our numerical
experiments, we observe that this choice leads to rather scattered
polynomial-degree distributions over the whole computational domain. The
choice~\eqref{eq_increm_p} is more conservative and leads to a smoother
overall polynomial-degree distribution. We believe that this choice is
preferable, at least as long as a polynomial-degree coarsening procedure is
not included in the adaptive loop. 
Another possibility is to use $\lceil\alpha p_{\ell,{\elm}}\rceil$ 
with $\alpha > 1$ instead of $p_{\ell,\elm} + \delta^\ver_\elm$, 
which corresponds to the theoretical developments 
in~\normalfont{\cite{Can_Noch_Stev_Ver_sat_16}}.
\er

\subsection{$hp$-decision on simplices} 

The second step in the module REFINE is to use the subsets $\Vmarkedh$
and $\Vmarkedp$ to decide whether $h$-, $p$- , or
$hp$-refinement should be performed on each simplex having at least one flagged vertex.
To this purpose, we define the following subsets:
\bse
\label{eq:decision_cell}\begin{align}
\Mellh &\eq \{ K\in \Th\;|\; \VK\cap \Vmarkedh \ne\emptyset \} \subset \Mell, \\
\Mellp &\eq \{ K\in \Th\;|\; \VK\cap \Vmarkedp \ne\emptyset \} \subset \Mell.
\end{align}
\ese
In other words, a simplex $\elm\in\Th$ is flagged for $h$-refinement
(resp., $p$-refinement) if it has at least one vertex flagged for
$h$-refinement (resp., $p$-refinement). Note that the subsets $\Mellh$ and
$\Mellp$ are not necessarily disjoint since a simplex can have some
vertices flagged for $h$-refinement and others flagged for $p$-refinement;
such simplices are then flagged for $hp$-refinement. Note also that
$\Mellh \cup \Mellp = \cup_{\ver \in\Vmarked} \Ta= \Mell$ is indeed
the set of marked simplices considered in the module MARK.

\subsection{$hp$-refinement}

In this last and final step, the subsets $\Mellh$ and $\Mellp$ are used
to produce first the next mesh $\Thh$ and then the next polynomial-degree
distribution $\tp_{\ell+1}$ on the mesh $\Thh$.

The next mesh $\Thh$ is a matching simplicial refinement of $\Th$ obtained by
dividing each flagged simplex $\elm \in \Mellh$ into at least two
simplices in a way that is consistent with the matching simplicial refinement
of $\Ta$ considered in Definition~\ref{def_h_ref} to build $\Tah$, i.e., such
that $\Tah \subset \Thh$ for all $\ver\in \Vmarkedh$. Note that to preserve
the conformity of the mesh, additional refinements beyond the set of flagged
simplices $\Mellh$ may be carried out when building $\Thh$. Several
algorithms can be considered to refine the mesh. In our numerical
experiments, we used the newest vertex bisection algorithm~\cite{Sew_NVB_72}.

After having constructed the next mesh $\Thh$, we assign the next
polynomial-degree distribution $\tp_{\ell+1}$ as follows. For all
$\widetilde{\elm}\in \Thh$, let $\elm$ denote its parent simplex in $\Th$. We
then set
\be \label{eq:new_p_unmarked} p_{\ell+1,\widetilde{\elm}} \eq
p_{\ell,\elm} \qquad \mbox{\textup{if}~$\elm\not\in\Mellp$},
\ee
that is, we assign the same polynomial degree to the children of a simplex
that is not flagged for $p$-refinement, whereas we set
\be \label{eq:new_p_marked} p_{\ell+1,\widetilde{\elm}} \eq
\max_{\ver\in \VK\cap \Vmarkedp} \big( p_{\ell,\elm} + \delta^\ver_\elm\big)
\qquad \mbox{\textup{if}~$\elm\in\Mellp$},
\ee
that is, we assign to the children of a simplex $\elm\in \Mellp$
flagged for $p$-refinement the largest of the polynomial degrees
considered in Definition~\ref{def_p_ref} to build the local residual liftings
associated with the vertices of $\elm$ flagged for $p$-refinement.

\subsection{Summary of the module REFINE}

The module REFINE is summarized in Algorithm~\ref{hp_strat_alg}.

\begin{algorithm}
\caption{(\textbf{module REFINE})}\label{hp_strat_alg}
\begin{algorithmic}[1]
\Module[{REFINE}]{$\Vmarked$}
\State \algorithmicrequire{ set of marked vertices $\Vmarked$}
\State \algorithmicensure{ next level mesh $\Thh$, polynomial-degree distribution $\tp_{\ell +1}$}
\ForAll{$\ver\in\Vmarked$}
\State Compute the $h$-refinement residual lifting $r^{\ver,h}$ from Definition~\ref{def_h_ref}
\State Compute the $p$-refinement residual lifting $r^{\ver,p}$ from Definition~\ref{def_p_ref}
\EndFor

\State $hp$-decision on vertices: build the subsets $\Vmarkedh$ and $\Vmarkedp$ from~\eqref{eq:decision_vertex}

\State $hp$-decision on simplices: build the subsets $\Mellh$ and $\Mellp$ from~\eqref{eq:decision_cell}

\State Build $\Thh$ from $\Th$ and $\Mellh$

\State Build $\tp_{\ell+1}$ on $\Thh$ from $\tp_\ell$,
$\{\delta^\ver_\elm\}_{\ver\in \Vmarkedp,\elm\in \Ta}$, and $\Mellp$
using~\eqref{eq:new_p_unmarked} and~\eqref{eq:new_p_marked}

\EndModule
\end{algorithmic}
\end{algorithm}

\begin{figure}[!h]
\begin{minipage}[b]{0.47\linewidth}
	\centering
	\includegraphics[width=\textwidth]{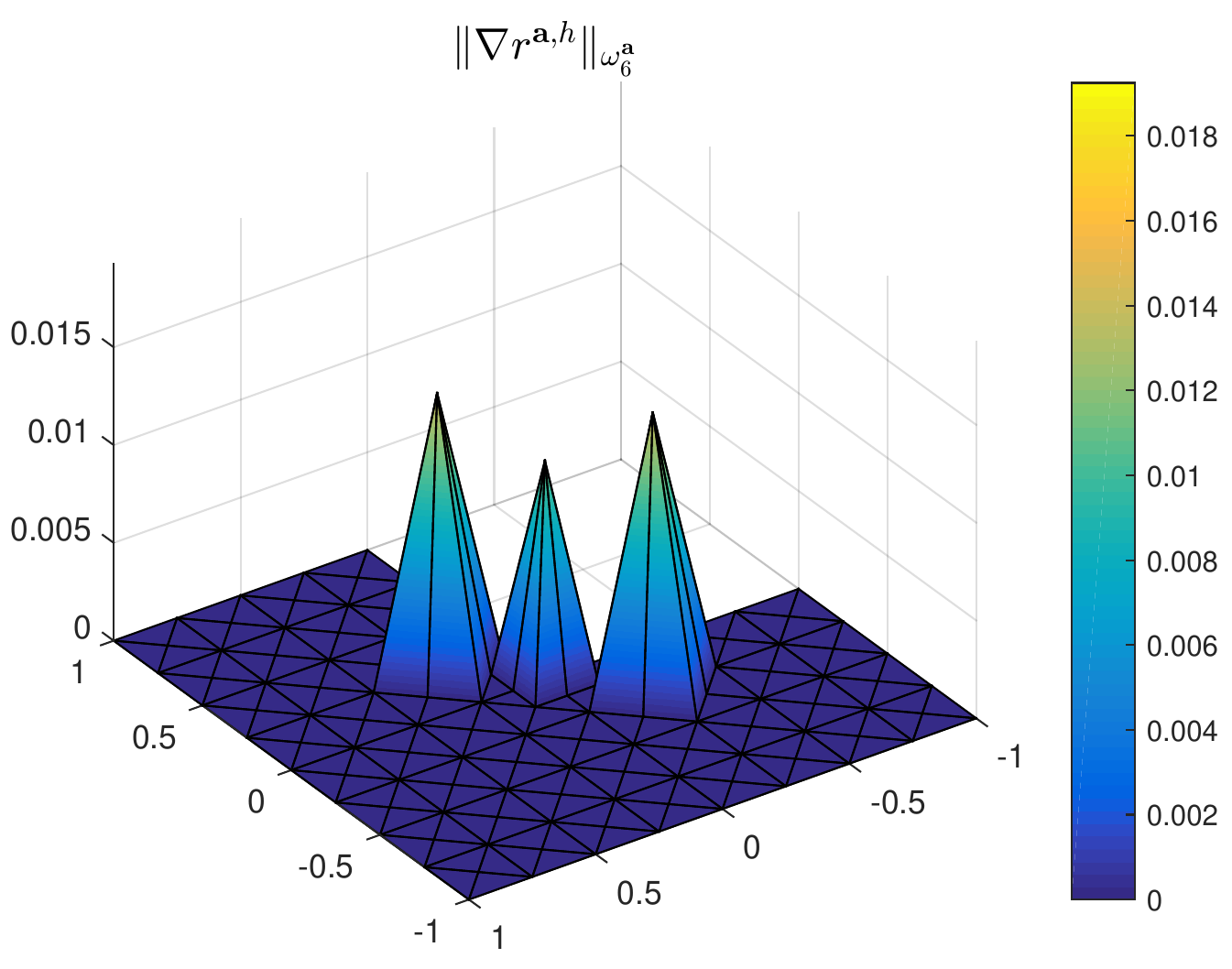}
\end{minipage}
\hfill
\begin{minipage}[b]{0.47\linewidth}
	\centering
	\includegraphics[width=\textwidth]{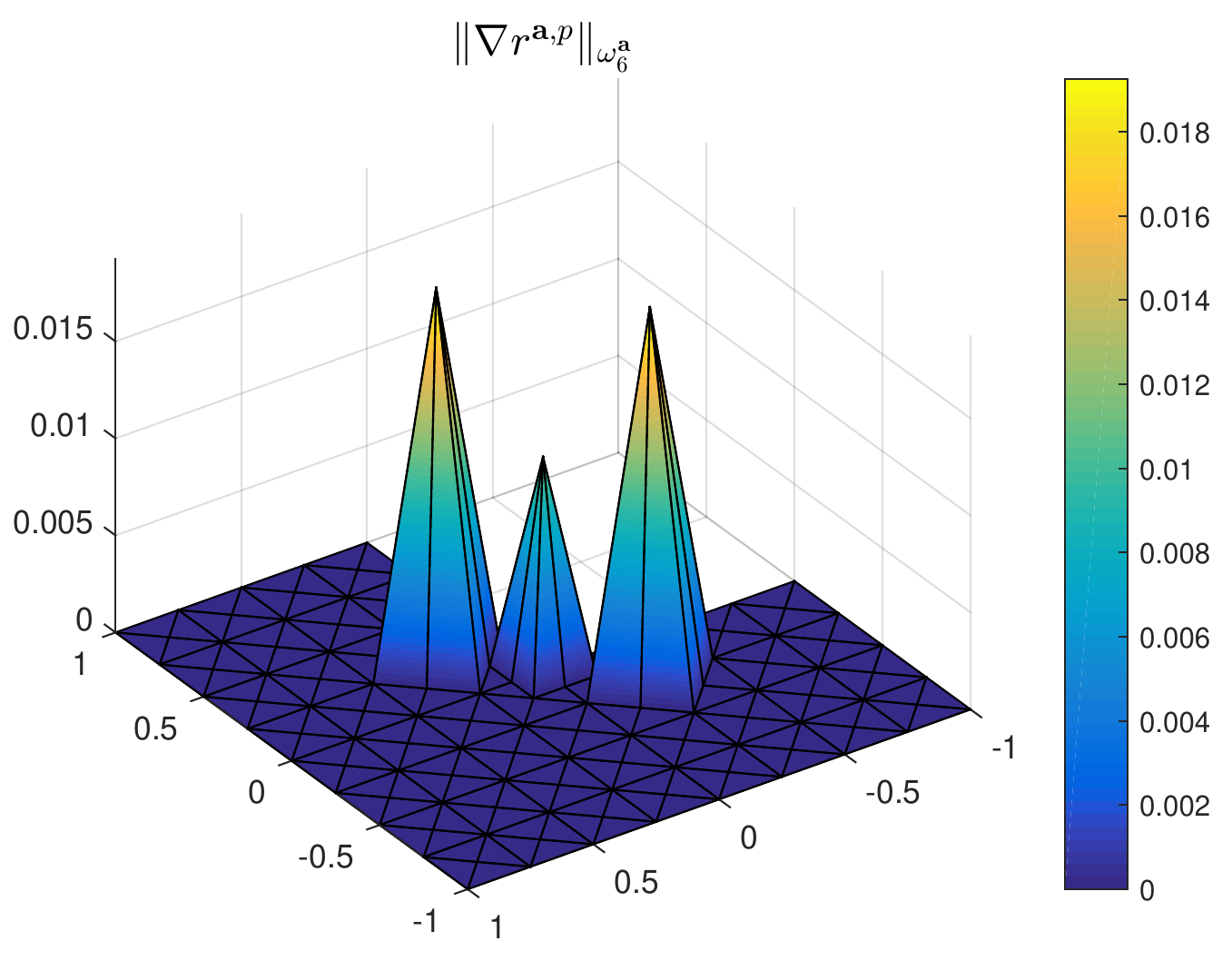}
\end{minipage}
	\caption{[L-shape problem from Section~\ref{sec_sing}] For the three marked vertices in \pd{$\widetilde{\V}_{6}^{\theta}$},
	we display the piecewise $\PP_1$ functions which take the value $\norm{\Gr r^{\ver,h}}_{\omega^\ver_{6}}$ in the vertex~$\ver$ and 0 elsewhere (\textit{left}) and the value
	$\norm{\Gr r^{\ver,p}}_{\omega^\ver_{6}}$ in the vertex~$\ver$ and 0 elsewhere (\textit{right}).}
	\label{fig:liftings_options}
\end{figure}

To illustrate Algorithm~\ref{hp_strat_alg}, we examine in detail a particular
situation with three marked vertices as encountered on the $6$th iteration
($\ell=6$) of the adaptive loop applied to the L-shape problem described in
Section~\ref{sec_sing} below. In
Figure~\ref{fig:chosen_lifting_decision} (left panel), we display the mesh
$\T_{6}$ and the polynomial-degree distribution $\tp_{6}$. There are three
marked vertices in \pd{$\widetilde{\V}_{6}^{\theta}$}. In
Figure~\ref{fig:liftings_options}, for the three marked vertices, we
visualize the norms $\norm{\Gr r^{\ver,h}}_{\omega^\ver_{6}}$ and $\norm{\Gr
r^{\ver,p}}_{\omega^\ver_{6}}$ which are the key ingredients for the
$hp$-decision on vertices. The resulting simplices flagged for $h$- and
$p$-refinement are shown in the central panel of
Figure~\ref{fig:chosen_lifting_decision}, whereas the right panel of
Figure~\ref{fig:chosen_lifting_decision} displays the next mesh $\T_{7}$ and
the next polynomial-degree distribution $\tp_{7}$.

\begin{figure}[!h]
\begin{minipage}[b]{0.3\linewidth}
	\centering
	\includegraphics[width=0.98\textwidth]{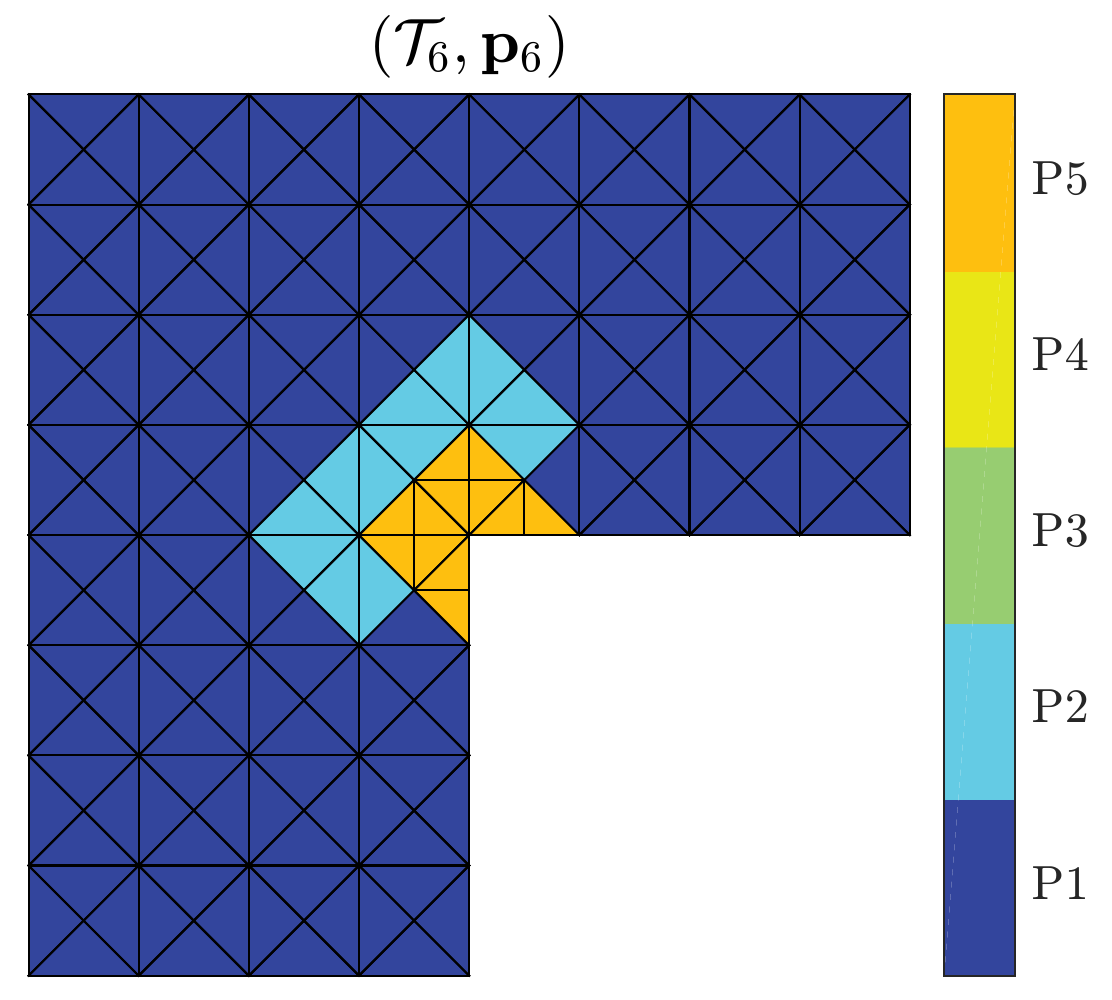}
\end{minipage}
\hfill
\begin{minipage}[b]{0.3\linewidth}
	\centering
	\includegraphics[width=1.05\textwidth]{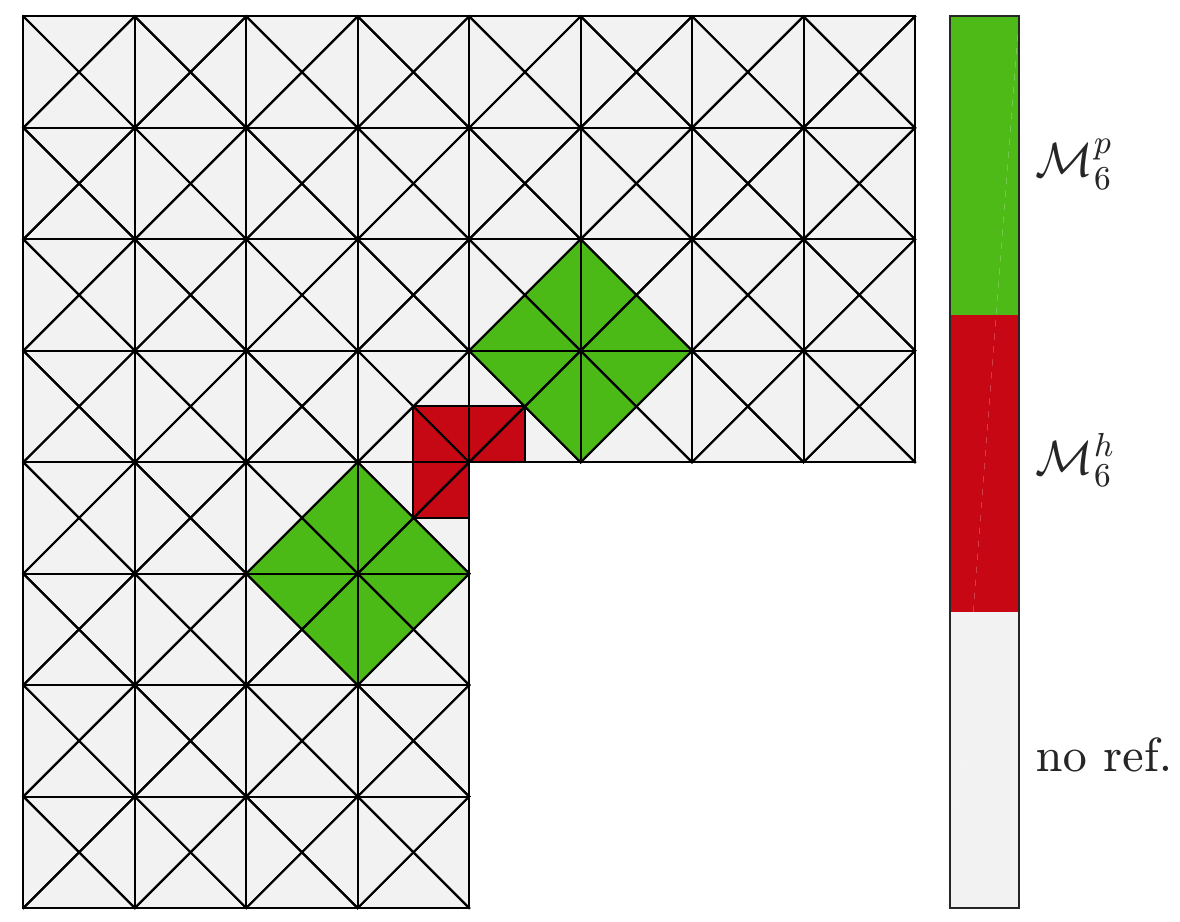}
\end{minipage}
\hfill
\begin{minipage}[b]{0.3\linewidth}
	\centering
	\includegraphics[width=0.98\textwidth]{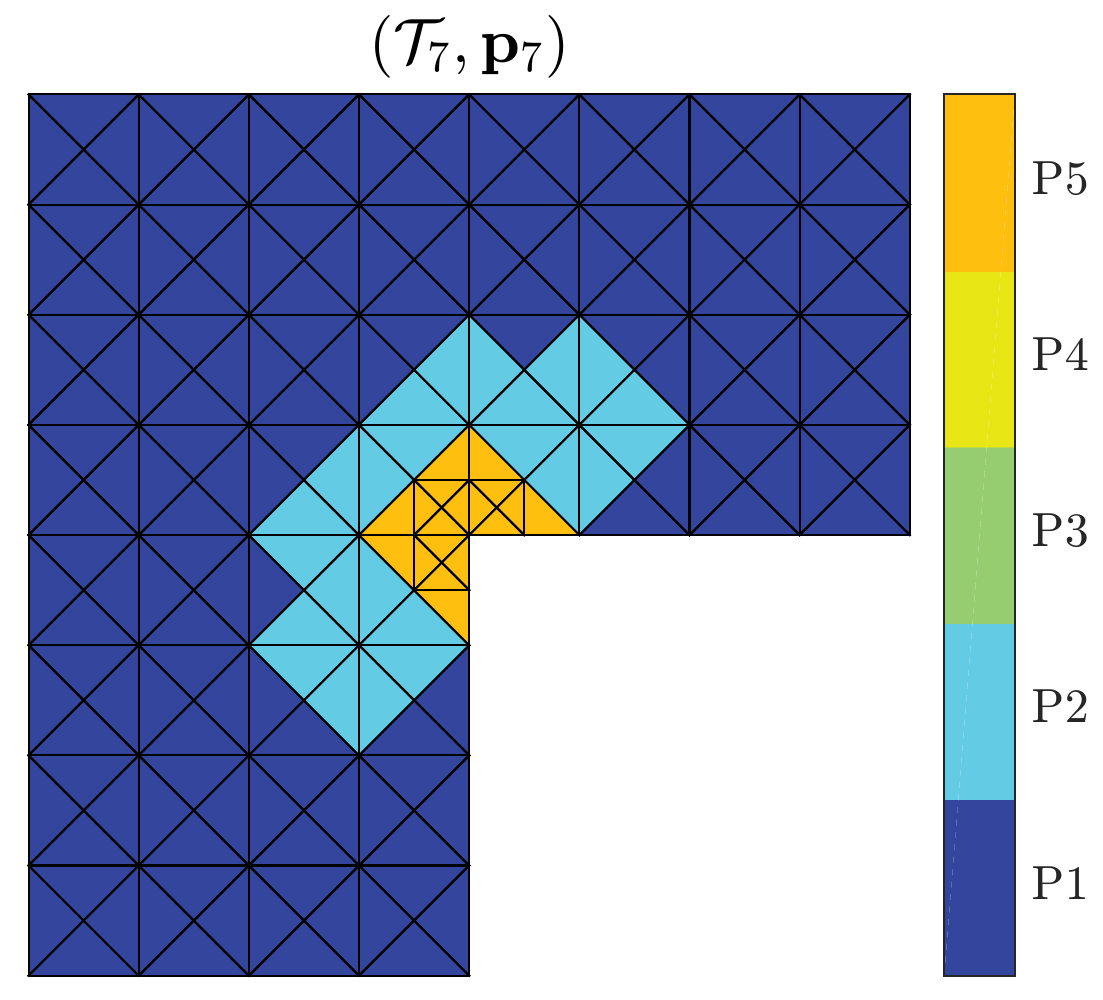}
\end{minipage}
	\caption{[L-shape problem from Section~\ref{sec_sing}]
		 The mesh and the polynomial degree distribution on the 6th iteration of the adaptive procedure~(\textit{left}). 	
		 Result of the $hp$-decision: simplices in $\M_{6}^h$ are shown in blue and simplices in $\M_{6}^p$ are shown in red, the two subsets $\M_{6}^h$ and $\M_{6}^p$ being here disjoint~(\textit{center}). The resulting mesh $\T_{7}$ and polynomial-degree distribution $\tp_{7}$~(\textit{right}).	
    }
	\label{fig:chosen_lifting_decision}
\end{figure}

\section{Guaranteed bound on the error reduction factor} \label{sec_err_red}

In this section we show that it is possible to compute, at marginal additional costs, a
guaranteed bound on the energy error reduction factor $\Cred$
from~\eqref{eq_C_red} on each iteration $\ell$
of the adaptive loop~\eqref{eq_afem_loop}. This bound can be computed right after the end of module
REFINE at the modest price of one additional primal solve in each
patch $\Ta$ associated with each marked vertex $\ver\in \Vmarked$.
Recall the set of marked simplices $\Mell = \cup_{\ver\in\Vmarked} \Ta$. 
Let us denote by $\omel \eq \cup_{\ver\in\Vmarked}\omega^\ver_\ell$ the corresponding
open subdomain; notice that a point $\tx$ is in $\overline{\omel}$ if and
only if there is $\elm\in \Mell$ so that $\tx \in \elm$. We start with
the following discrete lower bound result:

\bl[Guaranteed lower bound on the incremental error on marked simplices]
\label{lem_difference_bound_with_res} Let the mesh $\T_{\ell+1}$ and the
polynomial-degree distribution $\tp_{\ell+1}$ result from
Algorithm~\ref{hp_strat_alg}, and recall that $V_{\ell+1} =
\PP_{\tp_{\ell+1}}(\Thh) \cap \Hoo$ is the finite element space to be
used on iteration $(\ell+1)$ of the adaptive loop~\eqref{eq_afem_loop}. For
all the marked vertices $\ver \in \Vmarked$, let us set, in extension
of~\eqref{eq_V_a_h}, \eqref{eq_V_a_p},
\be \nonumber \Vahp \eq V_{\ell+1}|_{\omal}\cap
H^1_0(\omal), \ee and construct the residual lifting $\rahp \in \Vahp$ by
solving
\be \label{eq_local_fem_prbl_res}
(\Gr \rahp, \Gr v^{\ver,hp})_{\omal} = (f, v^{\ver,hp})_{\omal}-(\Gr u_{\ell}, \Gr v^{\ver,hp})_{\omal}\qquad \forall\, v^{\ver,hp} \in \Vahp.
\ee
Then, extending $\rahp$ by zero outside $\omal$, the following holds true: 	\be
\label{eq_res_norm_low_bound} \norm{\Gr(\uhh-\uh)}_\omel \geq \underline
\eta_{\Mell},\quad \underline \eta_{\Mell} \coloneqq \begin{cases}
\frac{\sum_{\ver \in \Vmarked}
    \norm[\big]{\Gr \rahp}_{\omal}^2 }{\norm[\Big]{\Gr\left(
    \sum_{\ver \in \Vmarked} \rahp\right)}_\omel} & \text{if } \sum_{\ver \in \Vmarked} \rahp\neq 0, \\ 0 & \text{otherwise.}
    \end{cases}
\ee \el

\bp Let $V_{\ell +1}(\omel)$ stand for the restriction of the space $V_{\ell
+1}$ to the subdomain $\omel$ and let
$V_{\ell +1}^{0}(\omel) \eq V_{\ell +1}(\omel)\cap H^1_0(\omel)$ stand for
the corresponding homogeneous Dirichlet subspace.
Note that $(\uhh-\uh)$ is a member of $V_{\ell +1}(\omel)$, but not necessarily of $V_{\ell +1}^{0}(\omel)$. Then, the following holds true:
\ban
    \norm{\Gr (\uhh - \uh)}_\omel & = \sup_{v_{\ell+1} \in V_{\ell
    +1}(\omel)} \frac{(\Gr (\uhh - \uh), \Gr v_{\ell+1})_\omel}{\norm{\Gr
    v_{\ell+1}}_\omel}\\
    & \geq \sup_{v_{\ell+1} \in V_{\ell +1}^{0}(\omel)} \dfrac{(\Gr (\uhh - \uh), \Gr v_{\ell+1})_\omel}{\norm{\Gr v_{\ell+1}}_\omel} \nonumber\\
    &= \sup_{v_{\ell+1} \in V_{\ell +1}^{0}(\omel)}\dfrac{(f,v_{\ell+1})_\omel-(\Gr
    \uh, \Gr v_{\ell+1})_\omel}{\norm{\Gr v_{\ell+1}}_\omel},
\ean
where we have used the definition~\eqref{eq_fem_model} of $\uhh$ on the mesh
$\T_{\ell+1}$, since $v_{\ell+1}$ extended by zero outside of $\omel$ belongs
to the space $V_{\ell + 1}$ whenever $v_{\ell+1} \in V_{\ell +1}^{0}(\omel)$. Now, choosing $v_{\ell+1} = \sum_{\ver \in \Vmarked} \rahp$
(note that this function indeed belongs to $V_{\ell +1}^{0}(\omel)$),
we infer that
\ifELSEVIER
\else
\ban
    \Bigg(f,\sum_{\ver \in \Vmarked} \rahp \Bigg)_\omel & - \Bigg(\Gr \uh, \Gr
    \Bigg( \sum_{\ver \in \Vmarked}\rahp \Bigg) \Bigg)_\omel \\
    &= \sum_{\ver \in \Vmarked} \big\{ \big(f, \rahp\big)_{\omal} - \big(\Gr \uh, \Gr \rahp \big)_{\omal}\big\}=\sum_{\ver \in \Vmarked}
    \norm[\big]{\Gr \rahp}_{\omal}^2,
\ean
%
\fi
where we have employed $\rahp$ as a test function
in~\eqref{eq_local_fem_prbl_res}. This finishes the proof. \ep

Our main result is summarized in the following contraction property in the
spirit of~\cite[Theorem~5.1]{Casc_Noch_cvg_non_res_11},
\cite[Proposition~4.1]{Can_Noch_Stev_Ver_sat_16}, and the references therein.
The specificity of the present work is that we obtain a guaranteed and computable bound on the error reduction factor. In contrast to these references, however, we
do not prove here that $\Cred$ is strictly smaller than one, although we
observe it numerically in Section~\ref{sec_num} below. 
\pd{We believe that one could show} $\Cred~<~1$ under
additional assumptions on the refinements, such as the interior node
property~\cite{Morin_converge_afem_02}, but we will not pursue this
consideration further here.

\bt[Guaranteed bound on the energy error reduction factor]
\label{th_contraction} Let the mesh $\T_{\ell+1}$ and the polynomial-degree
distribution $\tp_{\ell+1}$ result from Algorithm~\ref{hp_strat_alg}, and let
$V_{\ell+1} = \PP_{\tp_{\ell+1}}(\Thh) \cap \Hoo$ be the finite element space
to be used on iteration $(\ell+1)$ of the adaptive loop~\eqref{eq_afem_loop}.
Let $\underline \eta_{\Mell}$ be defined by~\eqref{eq_res_norm_low_bound}.
Then, unless $\eta(\Mell)=0$ in which case $\uh=u$ and the adaptive loop
terminates, the new numerical solution $\uhh \in V_{\ell + 1}$ satisfies
\be \label{eq_contr_property}
    \norm{\Gr (\uu -\uhh)} \le \Cred \norm{\Gr (\uu - \uh)} \quad \text{ with
    } \quad 0 \leq  \Cred \eq  \sqrt{1 - \theta^2 \frac{\underline
    \eta_{\Mell}^2}{\eta^2(\Mell)}} \leq 1.
\ee
\et

\bp We first observe that $\eta(\Mell)=0$ implies
using~\eqref{eq:def_Vmarked} and~\eqref{eq_est} that the error is zero
on iteration $\ell$, i.e., $\uu = \uh$, so that the adaptive
loop~\eqref{eq_afem_loop} terminates. Let us now assume that
$\eta(\Mell)\ne0$. Since the spaces $\{V_\ell\}_{\ell \geq 0}$ are
nested, cf.~\eqref{eq:nested}, Galerkin's orthogonality implies the following
Pythagorean identity:
\be  \nonumber
    \norm{\Gr (\uu - \uhh)}^2 = \norm{\Gr (\uu
    - \uh)}^2 - \norm{\Gr (\uhh - \uh)}^2.
\ee
Moreover, owing to Lemma~\ref{lem_difference_bound_with_res}, we infer
that
\[
    \norm{\Gr(\uhh-\uh)} \geq \norm{\Gr(\uhh-\uh)}_\omel \geq \underline \eta_{\Mell}
    = \frac{\underline \eta_{\Mell}}{\eta(\Mell)} \eta(\Mell).
\]
Using the marking criterion~\eqref{eq:def_Vmarked} and the definition of
$\Mell$, we next see that
\ifELSEVIER
\else
\ban 
	\norm{\Gr (\uu -\uhh)}^2 & \le \norm{\Gr
	    (\uu - \uh)}^2 - \frac{\underline \eta_{\Mell}^2}{\eta^2(\Mell)}
	    \eta^2(\Mell) \\
	    &\leq \norm{\Gr (\uu - \uh)}^2 - \theta^2
	    \frac{\underline \eta_{\Mell}^2}{\eta^2(\Mell)} \eta^2(\Th).
\ean
\fi
The assertion~\eqref{eq_contr_property} follows from the error estimate~\eqref{eq_est} and taking the square
root.
\ep

\begin{remark}[Local residual optimization] \label{rem_loc_res}
The use of the local residual liftings $r^{\ver,h}$ and $r^{\ver,p}$ from
Definitions~\ref{def_h_ref} and~\ref{def_p_ref} respectively in the
$hp$-decision criterion~\eqref{eq:decision_vertex} on marked vertices is motivated by
the result of Theorem~\ref{th_contraction}. Indeed, suppose that $r^{\ver,h}$
is larger than $r^{\ver,p}$ in norm, and that only $h$-refinement is
performed in the subdomain $\omal$ at the end of
Algorithm~\ref{hp_strat_alg}. Then, the local residual $r^{\ver,hp}$ from
Lemma~\ref{lem_difference_bound_with_res} coincides with $r^{\ver,h}$ which
means that by flagging the marked vertex $\ver$ for $h$-refinement, one
maximizes the contribution $\norm{\Gr \rahp}_{\omal}^2$ in the numerator
of~\eqref{eq_res_norm_low_bound} defining $\underline\eta_{\Mell}$.
It is also possible to design a more complex $hp$-refinement
strategy exploiting directly~\eqref{eq_contr_property}. Here we simply stick
to Algorithm~\ref{hp_strat_alg} which in our numerical experiments reported
in Section~\ref{sec_num} below leads to exponential convergence rates.
\end{remark}

\begin{remark}[A sharper bound] \label{rem_sharp}
Theorem~\ref{th_contraction} obviously also holds true with the slightly sharper
constant $\Cred=\sqrt{1-\frac{\underline
\eta_{\Mell}^2}{\eta^2(\T_{\ell})}}$. This is equivalent to considering
in~\eqref{eq_contr_property} $\theta_\ell$ such that $\eta(\Mell) =
\theta_\ell \,\eta({\Th})$ in place of $\theta$, a strategy adopted in the
numerical experiments in Section~\ref{sec_num} below.
We note that $\theta_{\ell} \ge \theta$, however employing $\theta_{\ell}$ in Algorithm~\ref{def_mark_patch_alg} would lead to the same set of marked simplices $\Mellth~=~\Mell$.
\end{remark}

\section{Numerical experiments} \label{sec_num}

We consider two test cases for the model problem~\eqref{eq_weak_form}, both in two space dimensions, one with a (relatively) smooth weak solution and one with a singular weak solution. Our main goal with the numerical experiments is to verify that the $hp$-refinement strategy of Algorithm~\ref{hp_strat_alg}
leads to an exponential rate of convergence with respect to the number of
degrees of freedom $\mathrm{DoF}_\ell$ of the finite element spaces $V_\ell$
in the form
\be
\label{eq_exp_converg}
\norm{\Gr (\uu - \uh)} \le C_1 \exp\left( {-C_2 \mathrm{DoF}_\ell^{\frac 1 3}}\right),
\ee
with positive constants $C_1$, $C_2$ independent of $\mathrm{DoF}_\ell$. In
addition, we assess the sharpness of the guaranteed bound on the
reduction factor $\Cred$ from Theorem~\ref{th_contraction} by means of the
effectivity index defined as
\be \label{eq_eff_index}
    I_{\mathrm{red}}^{\mathrm{eff}} =
    \frac{\Cred}{\frac{\norm{\Gr(\uu-\uhh)}}{\norm{\Gr(\uu-\uh)}}}.
\ee
We always consider the (well-established) choice $\theta = 0.5$ for the
marking parameter, fine-tuning it on each step to $\theta_\ell$ as
described in Remark~\ref{rem_sharp}. As mentioned above, we apply the newest
vertex bisection algorithm~\cite{Sew_NVB_72} to perform $h$-refinement and we
use the polynomial-degree increment~\eqref{eq_increm_p} to perform
$p$-refinement.

We compare the performance of our $hp$-refinement algorithm to two other
algorithms based on a different $hp$-decision criteria, namely the
\textsf{PARAM} and \textsf{PRIOR} criteria from the survey
paper~\cite{Mitchell_hp_review_14} which are both based on a local smoothness
estimation. These criteria hinge on the local $L^2$-orthogonal projection
$\uhp$ of the numerical solution $\uh$ onto the local lower-polynomial-degree
space $\PP_{p_{\ell,\elm} - 1}(\elm)$ for all the marked simplices
$K\in\Mell$. This leads to the local quantity $\eta_{\elm}^{p-1} \eq
\norm{\Gr (\uh - \uhp)}_\elm$; in case of $p_{\ell,\elm} = 1$, when the quantity~$\eta_{\elm}^{p-1}$ is not available, for both criteria, the marked simplex $\elm$ is $p$-refined. The criterion
\textsf{PARAM}~\cite{Gui_Bab_hp_III_86} relies on the local smoothness
indicator $g_\elm \eq \eta_\elm/\eta_{\elm}^{p-1}$ and a user-defined
parameter $\gamma>0$; the marked simplex $\elm$ is $h$-refined if $g_\elm >
\gamma$, and otherwise it is $p$-refined. The presence of the parameter
$\gamma$ is a drawback of this criterion; in our experiments we use the
values $\gamma = 0.3$ and $\gamma = 0.6$, as suggested
in~\cite{Mitchell_hp_review_14}. The criterion \textsf{PRIOR}, which is a
simplified version of the one proposed in~\cite{Sul_Hous_Schw_hp_hyper_00},
relies on the quantity $s_\elm \eq 1 - \log({\eta_\elm/\eta_\elm^{p-1}}) /
\log(p_{\ell,\elm}/(p_{\ell,\elm} -1))$; the marked simplex $\elm$ is $h$-refined if
$p_{\ell,\elm} > s_\elm - 1$, and otherwise it is $p$-refined. To make the
comparison with our approach more objective, we apply for both criteria the
suggested $p$-refinement only to those simplices such that $p_{\ell,\elm} =
\min_{\elm' \in \T_\ver} p_{\ell,\elm'}$.

\subsection{Smooth solution (sharp Gaussian)} \label{sec_reg}

We consider a square domain $\Om = (-1,1) \times (-1,1)$ and a weak solution
that is smooth but has a rather sharp peak
\be \nonumber
    \uu (x,y) = (x^2-1)(y^2-1)\exp{(-100 (x^2+y^2))}.
\ee
We start from a criss-cross initial mesh $\T_0$ with $\max_{\elm \in
\T_0}h_\elm = 0.25$ and a uniform polynomial-degree distribution equal to $1$ on all triangles.

\begin{figure}[!ht]
\begin{minipage}[b]{0.49\linewidth}
	\centering
	\includegraphics[width=0.77\textwidth]{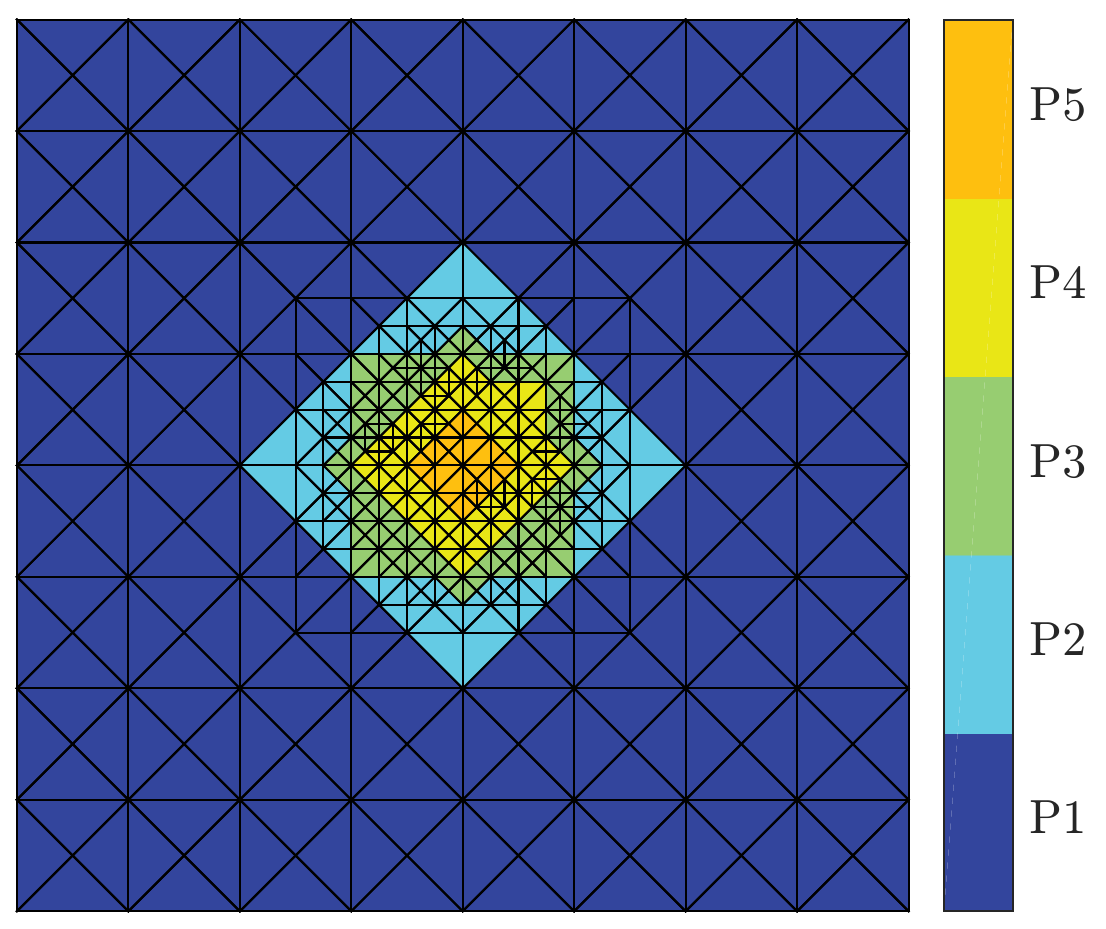}
\end{minipage}
\begin{minipage}[b]{0.49\linewidth}
	\centering
	\includegraphics[width=0.85\textwidth]{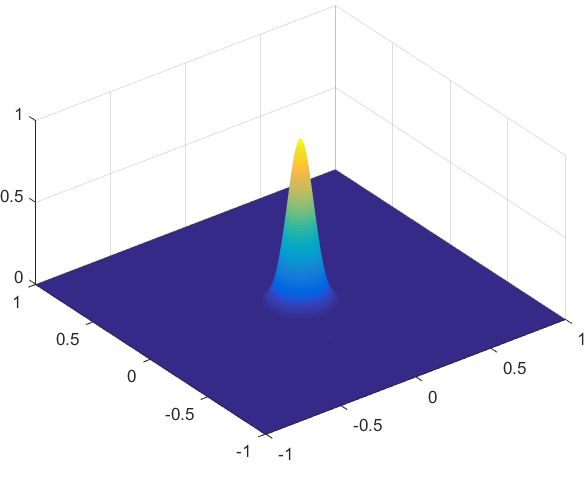}
\end{minipage}
	\caption{[{Sharp-Gaussian} of Section~\ref{sec_reg}]
    The final mesh and polynomial-degree distribution
     obtained after $30$ iterations of the $hp$-adaptive procedure (\textit{left}) and the obtained numerical solution $\uu_{30}$ (\textit{right}).}
	\label{fig:peak_meshes}
\end{figure}

\begin{figure}[!ht]
	\centering
	\includegraphics[width=0.40\textwidth]{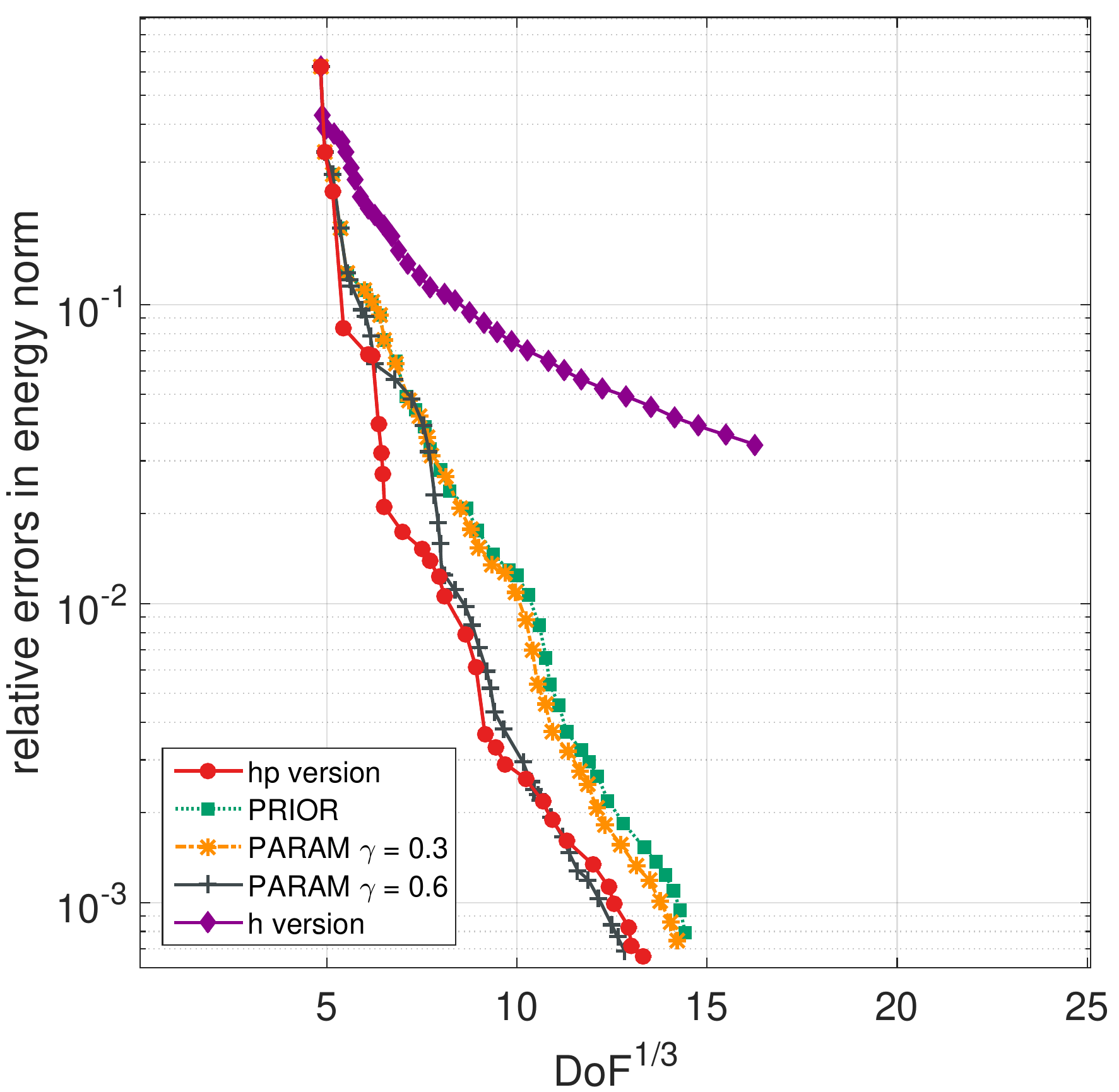}
    \caption{\pd{[{Sharp-Gaussian} of Section~\ref{sec_reg}]
            Relative energy error $\norm{\Gr (\uu - \uh)}/\norm{\Gr \uu}$ as a function
            of $\mathrm{DoF}_\ell^{\frac 1 3}$, obtained using the present $hp$-decision criterion, the criteria \textsf{PRIOR} and \textsf{PARAM} ($\gamma = 0.3$, $\gamma = 0.6$), and using only $h$-refinement.}}
	\label{fig:peak_conv}
\end{figure}

\begin{center}
\begin{figure}[!ht]
	\centering
	\includegraphics[width=0.4\textwidth]{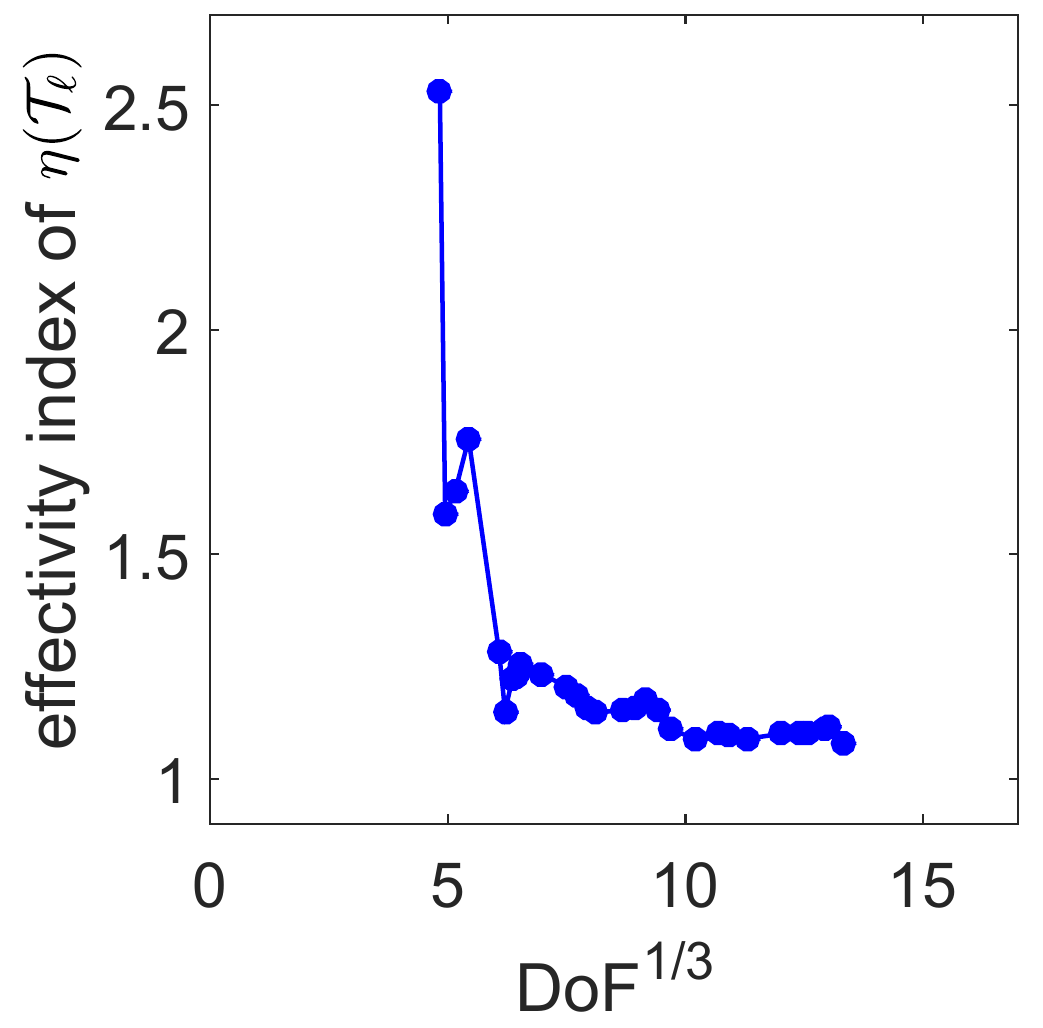}
	\caption{\pd{[{Sharp-Gaussian} of Section~\ref{sec_reg}] Effectivity indices of the error estimators $\eta(\Th)$ from Theorem~\ref{thm_est}, defined as the ratio $\eta(\Th)/\norm{\Gr (\uu - \uh)}$, throughout the $hp$-adaptive procedure.}
    }
	\label{fig:peak_Ieff_estim}
\end{figure}
\end{center}

\begin{center}
\begin{figure}[!ht]
\begin{minipage}[b]{0.49\linewidth}
	\hfill
	\includegraphics[width=0.81\textwidth]{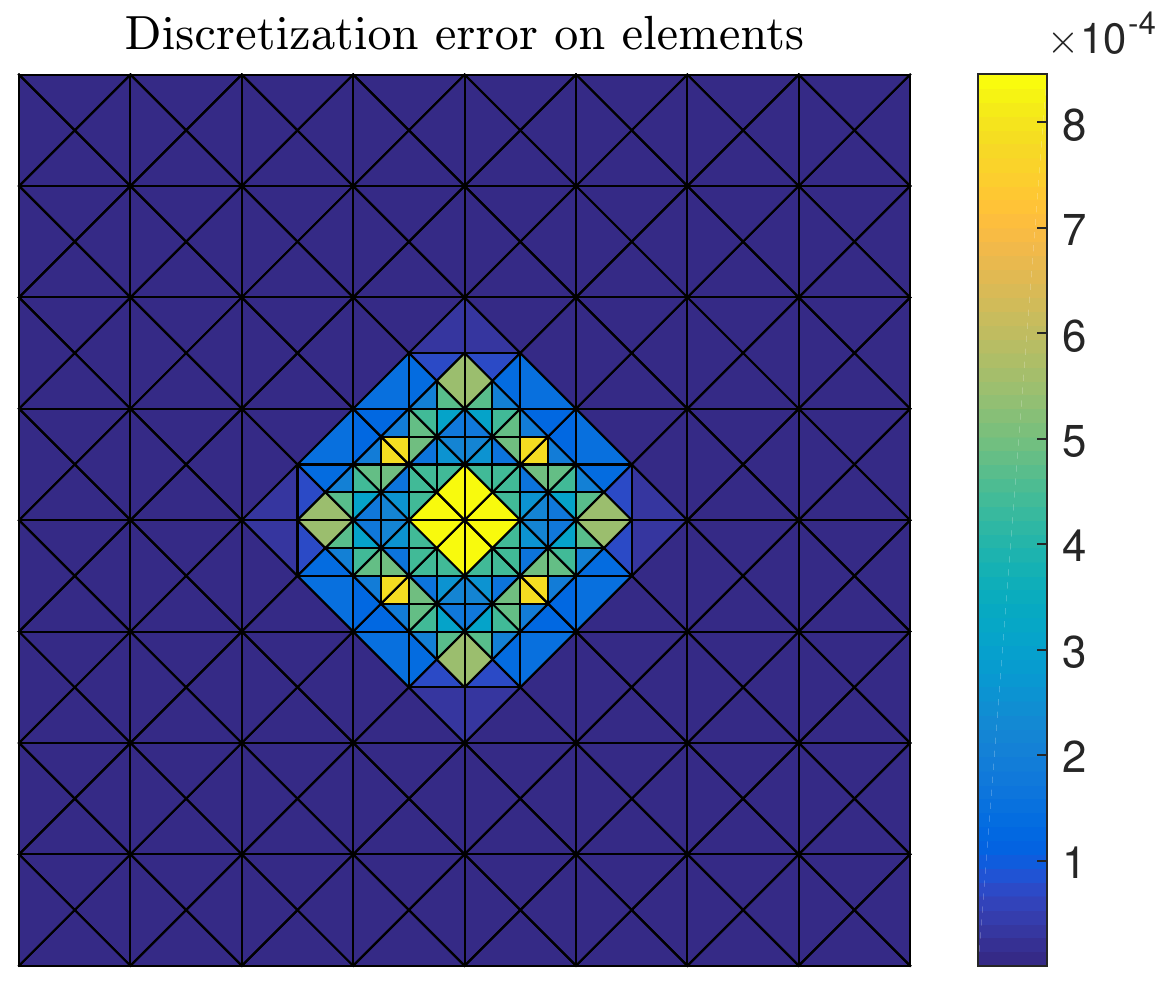}
\end{minipage}
\hfill
\begin{minipage}[b]{0.49\linewidth}
	\centering
	\includegraphics[width=0.80\textwidth]{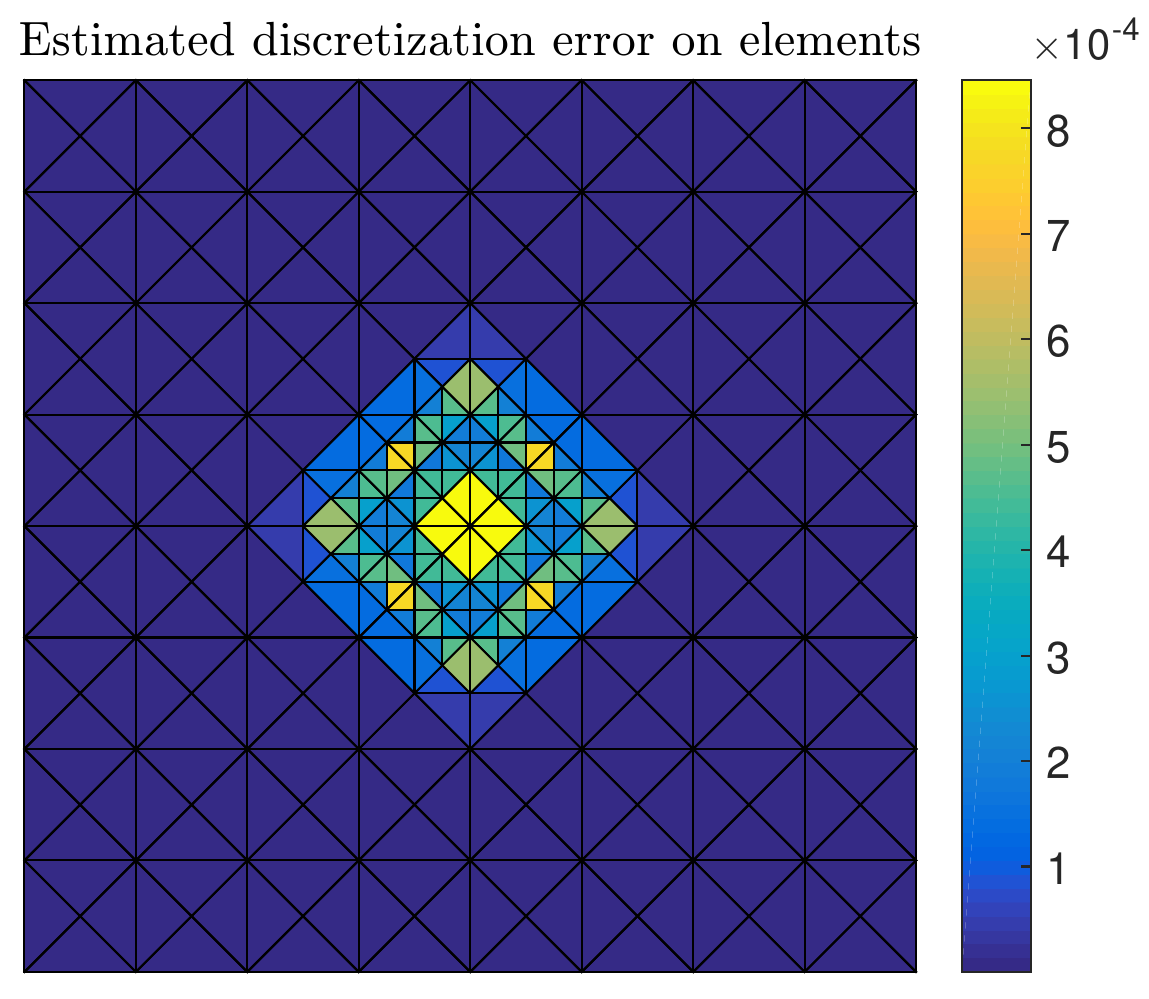}
\end{minipage}
	\caption{[{Sharp-Gaussian} of Section~\ref{sec_reg}] The distribution of the
	energy error $\norm{\Gr(\uu-\uh)}_\elm$~(\textit{left}) and of the error estimators
	$\eta_\elm$ from Theorem~\ref{thm_est}~(\textit{right}), $\ell = 20$. The effectivity index of the estimate defined as $\eta({\T_{20}}) / \norm{\Gr(\uu - \uu_{20})}$ is $1.1108$.
    }
	\label{fig:peak_estim}
\end{figure}
\end{center}

\begin{figure}[!ht]
\begin{minipage}[b]{0.49\linewidth}
	\centering
	\includegraphics[width=0.85\textwidth]{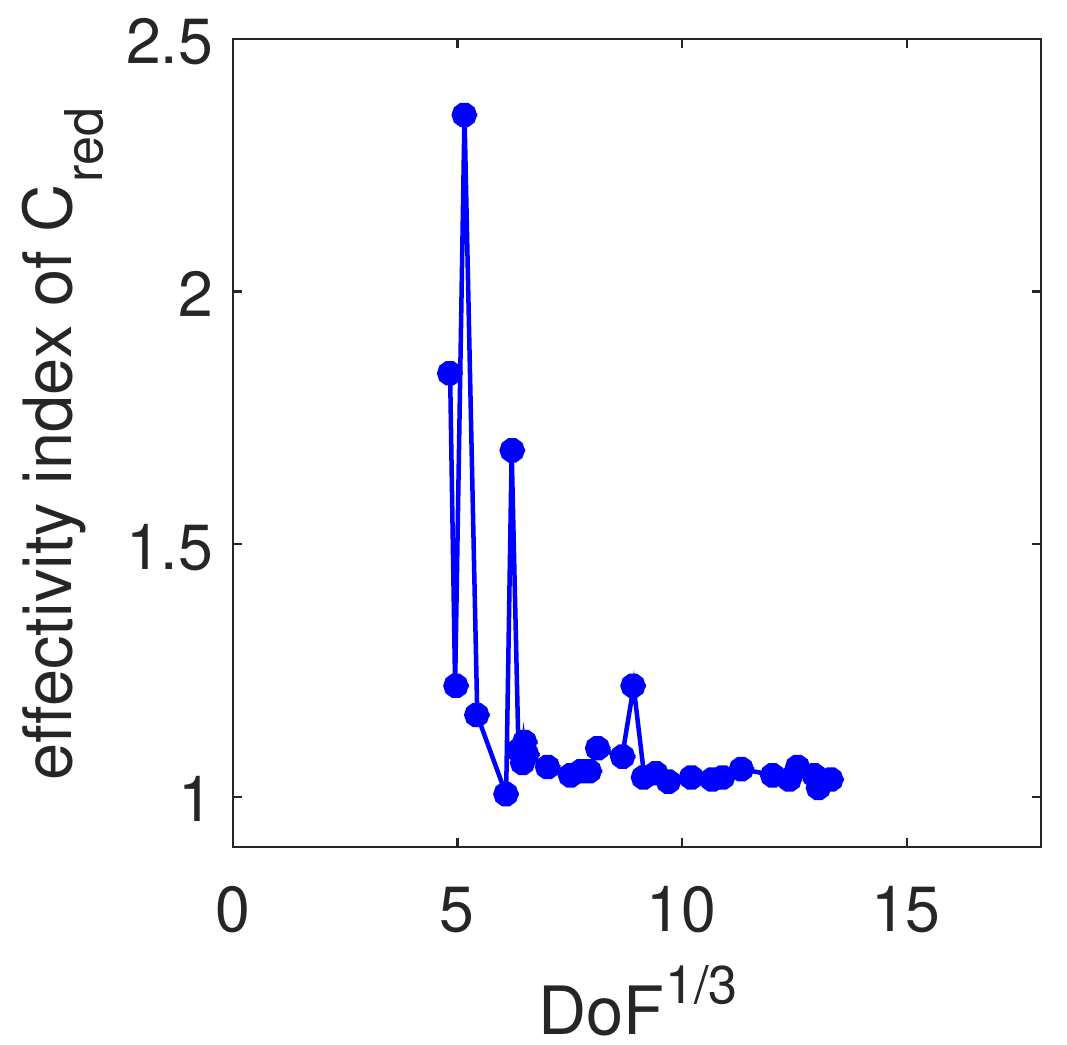}
\end{minipage}
\hfill
\begin{minipage}[b]{0.49\linewidth}
	\centering
	\includegraphics[width=0.8\textwidth]{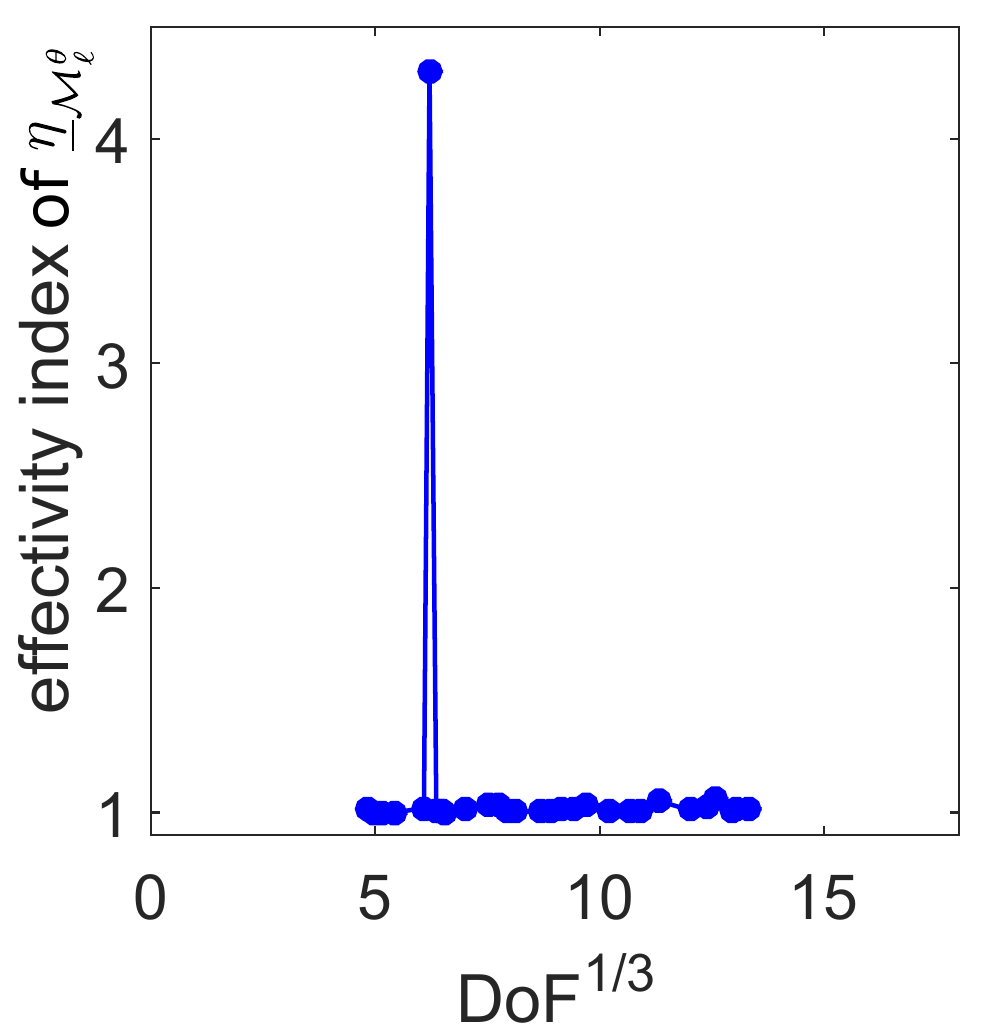}
\end{minipage}
	\caption{[{Sharp-Gaussian} of Section~\ref{sec_reg}] Effectivity indices~\eqref{eq_eff_index}
    for the error reduction factor $\Cred$ from Theorem~\ref{th_contraction} (\textit{left}) and effectivity indices
    for the lower bound $\underline \eta_{\Mell}$ from
	Lemma~\ref{lem_difference_bound_with_res} defined as the ratio
	$\norm{\Gr(\uhh-\uh)}_\omel/\underline \eta_{\Mell}$
    (\textit{right}).}
	\label{fig:peak_Ieff_contr_low_bound}
\end{figure}

Figure~\ref{fig:peak_meshes} presents the final mesh and polynomial-degree distribution obtained after 30 steps of the $hp$-adaptive procedure~\eqref{eq_afem_loop} (left panel) along with the obtained numerical solution (right panel).
Figure~\ref{fig:peak_conv} displays the relative error $\norm{\Gr(\uu-\uh)}/\norm{\Gr \uu}$ as a function of $\mathrm{DoF}_\ell^{\frac 1 3}$ in logarithmic-linear scale to illustrate that the present $hp$-adaptive procedure leads to an asymptotic exponential rate of convergence.
The values of \pd{the} constants $C_{1}$ and $C_{2}$ from~\eqref{eq_exp_converg} given by the 2-parameter least squares fit are \num[round-mode=places,round-precision=2]{3.96636} and \num[round-mode=places,round-precision=2]{0.695854}, respectively.
The value of $C_{2}$ indicates the slope steepness of the fitted line in logarithmic-linear scale, in particular, the higher value of $C_{2}$, the steeper slope.
For comparison, we also plot the relative error obtained when using the $hp$-decision criteria \textsf{PRIOR} and \textsf{PARAM} described above and also for the pure $h$-version of the adaptive loop.
\pd{The quality of the a posteriori error estimators of Theorem~\ref{thm_est} throughout the whole $hp$-adaptive process can be appreciated in Figure~\ref{fig:peak_Ieff_estim} where the effectivity indices, defined as the ratio of the error estimator~$\eta(\Th)$ and the actual error~$\norm{\Gr (\uu - \uh)}$, are presented. 
  Then, in Figure~\ref{fig:peak_estim} we compare the actual and estimated error distributions on iteration $\ell=20$ of the adaptive loop, showing excellent agreement.}  
Figure~\ref{fig:peak_Ieff_contr_low_bound} (left panel) presents the effectivity index for the reduction factor $\Cred$, see~\eqref{eq_eff_index}, throughout the adaptive process.
Overall, values quite close to one are obtained, except at some of the first iterations where the values are larger but do not exceed $2.5$.
Moreover, all the values are larger than one, confirming that the bound on the reduction factor $\Cred$ is indeed guaranteed.
Figure~\ref{fig:peak_Ieff_contr_low_bound} (right panel) examines the quality of the lower bound $\underline \eta_{\Mell}$ from Lemma~\ref{lem_difference_bound_with_res} by plotting the ratio of the left-hand side to the right-hand side of the lower bound in~\eqref{eq_res_norm_low_bound}.
Except for one iteration where this ratio takes a larger value close to $4.5$, we observe that this ratio takes always values quite close to, and larger than, one, indicating that $\underline \eta_{\Mell}$ delivers a sharp and guaranteed lower bound on the energy error decrease.
To give some further insight into the proposed $hp$-adaptive process, we present in Tables~\ref{table:history_peak_1} and~\ref{table:history_peak_2} some details on the $hp$-refinement decisions throughout the first 10 and the last 10 iterations of the adaptive loop.
Finally, Table~\ref{table:comp_iter_dofs} (top) compares the different
strategies namely in terms of the number of iterations of the adaptive
loop~\eqref{eq_afem_loop}; here our strategy is a clear winner.

\begin{table}[h!]
\centering
\begin{tabular}{l cccccccccc }
 \toprule
 Iteration& 1 & 2 & 3 & 4 & 5 & 6 & 7 & 8 & 9 & 10 \\ \hline
 Triangles  & 256 & 256 & 256 & 256 & 264 & 264 & 264 & 264 & 264 & 264 \\
 Maximal polynomial degree  & 1 & 2 & 3 & 4 & 4 & 4 & 4 & 4 & 4 & 4  \\
 Marked vertices   & 1 & 1 & 1 & 1 & 2 & 2 & 2 & 1 & 1 & 1  \\
 Triangles flagged for $h$-refinement  & 0 & 0 & 0 & 8 & 0 & 0 & 0 & 0 & 0 & 8 \\
 Triangles flagged for $p$-refinement & 8 & 8 & 8 & 0 & 12 & 12 & 4 & 2 & 2 & 0  \\
 Triangles flagged for $hp$-refinement & 0 & 0 & 0 & 0 & 0 & 0 & 0 & 0 & 0 & 0 \\
 \hline
 \bottomrule
\end{tabular}
\caption{[{Sharp-Gaussian} of Section~\ref{sec_reg}] Refinement decisions in Algorithm~\ref{hp_strat_alg} during the first 10 iterations of the adaptive loop~\eqref{eq_afem_loop}.
}
\label{table:history_peak_1}
\end{table}

\begin{table}[h!]
\centering
\begin{tabular}{l cccccccccc }
 \toprule
 Iteration &  20 & 21	& 22	& 23	& 24	& 25	& 26	& 27	& 28	& 29	\\ \hline
 Triangles  &  392 & 406	& 430	& 450	& 478	& 514	& 552	& 580	& 612	& 612	\\
 Maximal polynomial degree  &  5 & 5	& 5	& 5	& 5	& 5	& 5	& 5	& 5	& 5	\\
Marked vertices   &  2 & 3	& 3	& 3	& 4	& 3	& 2	& 3	& 4	& 4	\\
 Triangles flagged for $h$-refinement  &  12 & 24	& 16	& 24	& 30	& 23	& 14	& 21	& 0	& 28	\\
 Triangles flagged for $p$-refinement &  0 & 0	& 4	& 0	& 0	& 0	& 0	& 0	& 16	& 0	\\
 Triangles flagged for $hp$-refinement &  0 & 0	& 0	& 0	& 0	& 0	& 0	& 0	& 0	& 0\\
 \hline
 \bottomrule
\end{tabular}
\caption{[{Sharp-Gaussian} of Section~\ref{sec_reg}] Refinement decisions in Algorithm~\ref{hp_strat_alg} during the last 10 iterations of the adaptive loop~\eqref{eq_afem_loop}.
}
\label{table:history_peak_2}
\end{table}

\FloatBarrier

\subsection{Singular solution (L-shape domain)} \label{sec_sing}

In our second test case, we consider the L-shape domain $\Om = (-1,1)
\times(-1,1)\, \setminus\, [0,1] \times [-1,0]$ with $f = 0$ and
the exact solution (in polar coordinates)
\[
    u(r,\varphi) = r^{\frac{2}{3}} \sin\left(\frac{2 \varphi}{3}\right).
\]
\pd{For this test case, following~\cite[Theorem~3.3]{hp_Dol_Ern_Voh_16} and the references therein, 
the error estimator $\eta(\Th)$ employed within the adaptive procedure takes into account also 
the error from the approximation of the inhomogeneous Dirichlet boundary condition prescribed by the exact solution on $\pt \Om$.}
We start the computation on a criss-cross grid \pd{$\T_0$ with $\max_{\elm \in
\T_0}h_\elm = 0.25$} and all the polynomial degrees set uniformly to $1$.

\begin{figure}[!ht]
\centering
\begin{minipage}[b]{0.41\linewidth}
	\centering
	\includegraphics[width=0.9\textwidth]{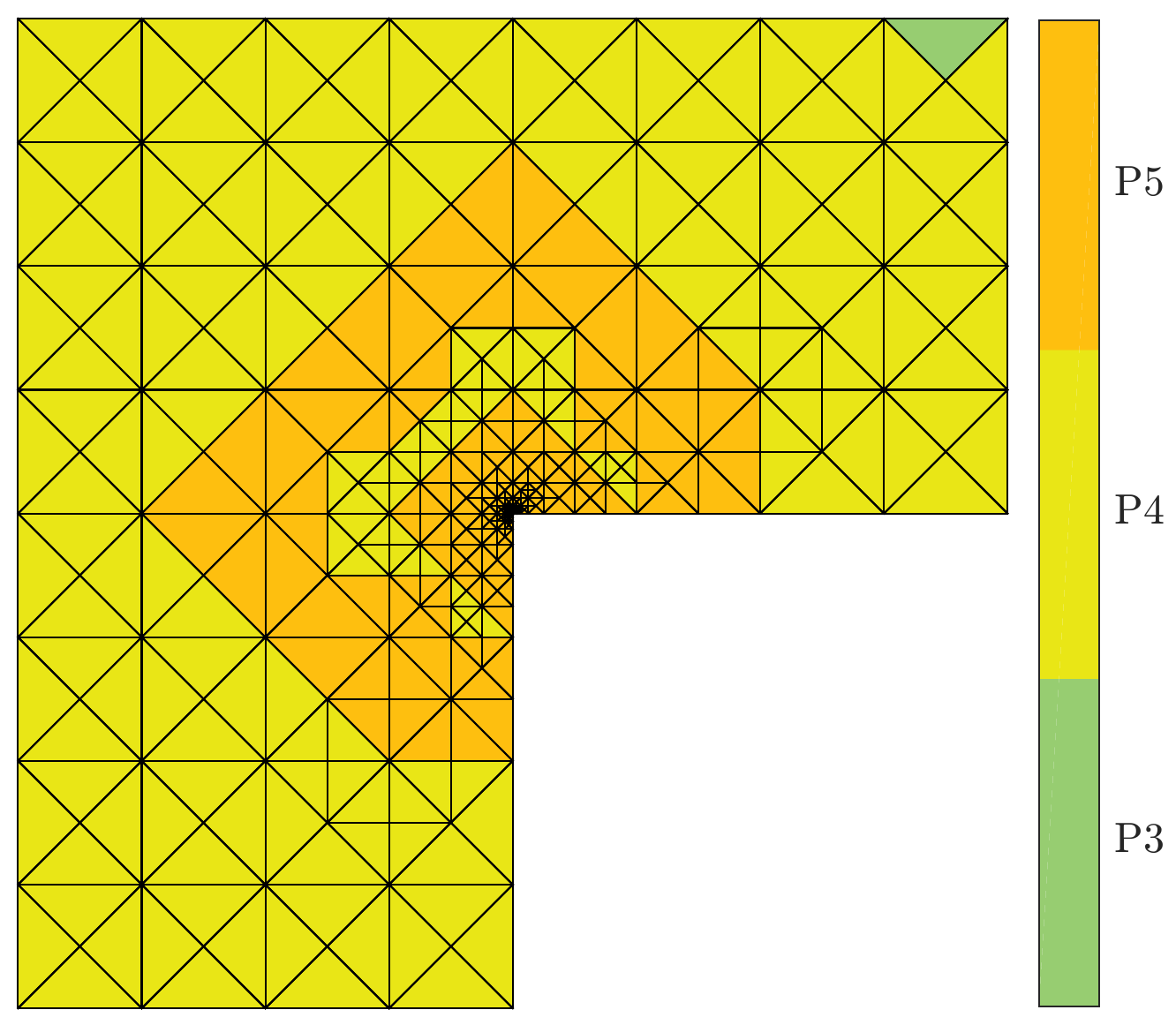}
\end{minipage}
\hfill
\begin{minipage}[b]{0.41\linewidth}
	\centering
	\includegraphics[width=0.82\textwidth]{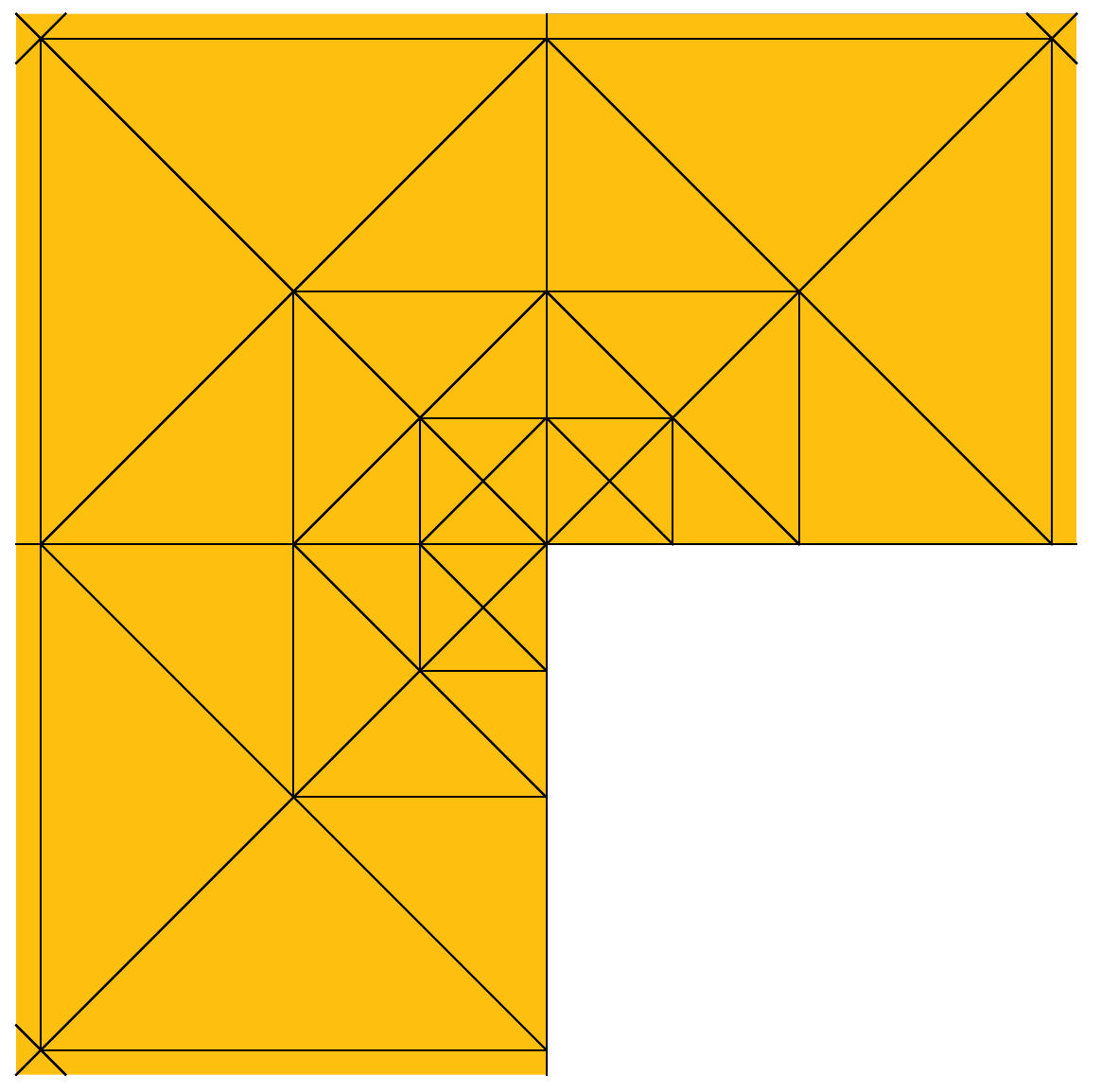}
\end{minipage}
	\caption{[L-shape domain of Section~\ref{sec_sing}] The final mesh and polynomial-degree distribution obtained after $65$ iterations of the $hp$-adaptive procedure (\textit{left}) and a zoom in $[-10^{-6},10^{-6}] \times [-10^{-6},10^{-6}]$ near the re-entrant corner (\textit{right}).
    }
	\label{fig:Final mesh}
\end{figure}

\begin{figure}[!ht]
\begin{minipage}[b]{0.49\linewidth}
	\centering
	\includegraphics[width=0.83\textwidth]{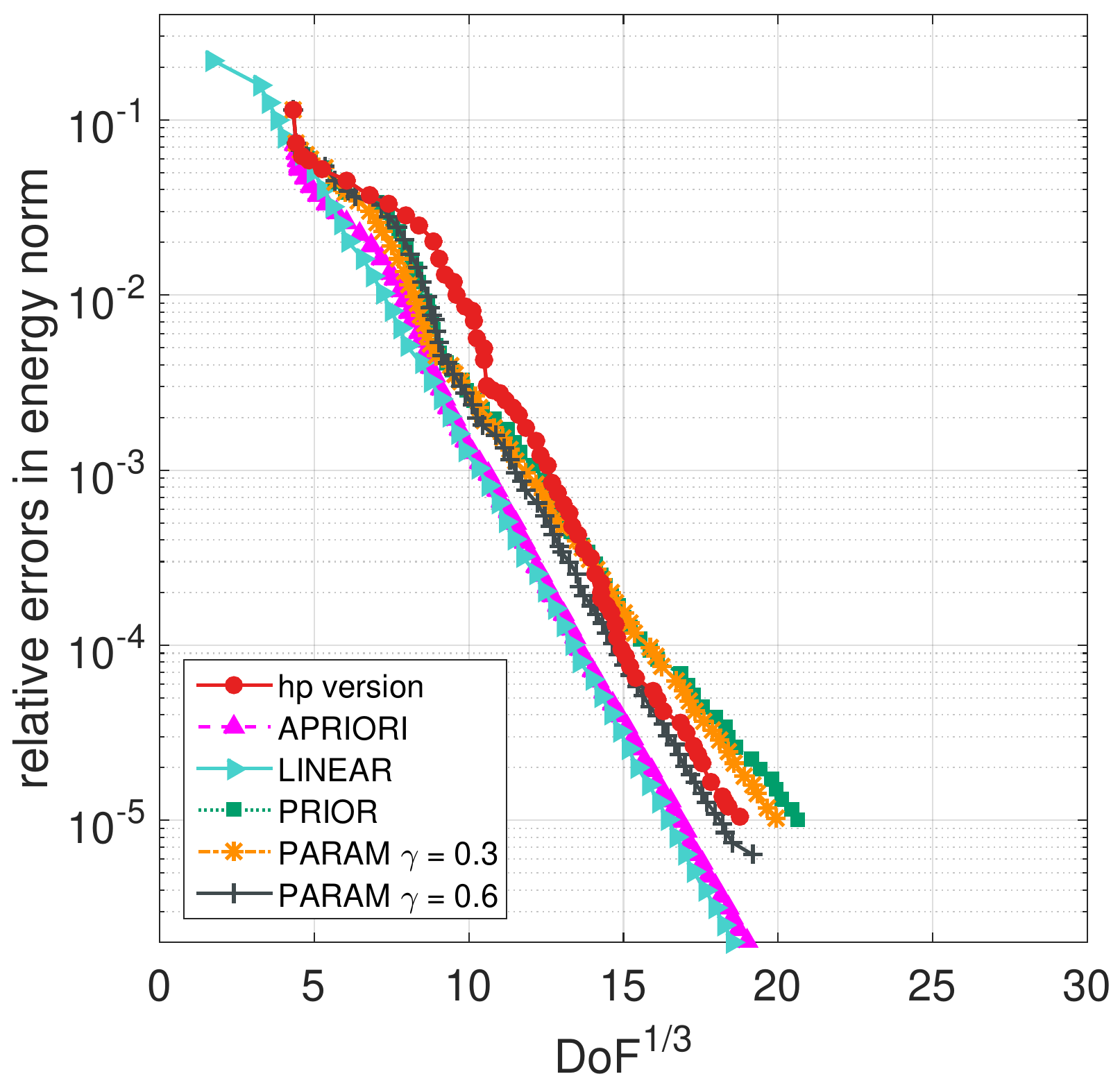}
\end{minipage}
\hfill
\begin{minipage}[b]{0.49\linewidth}
	\centering
	\includegraphics[width=0.83\textwidth]{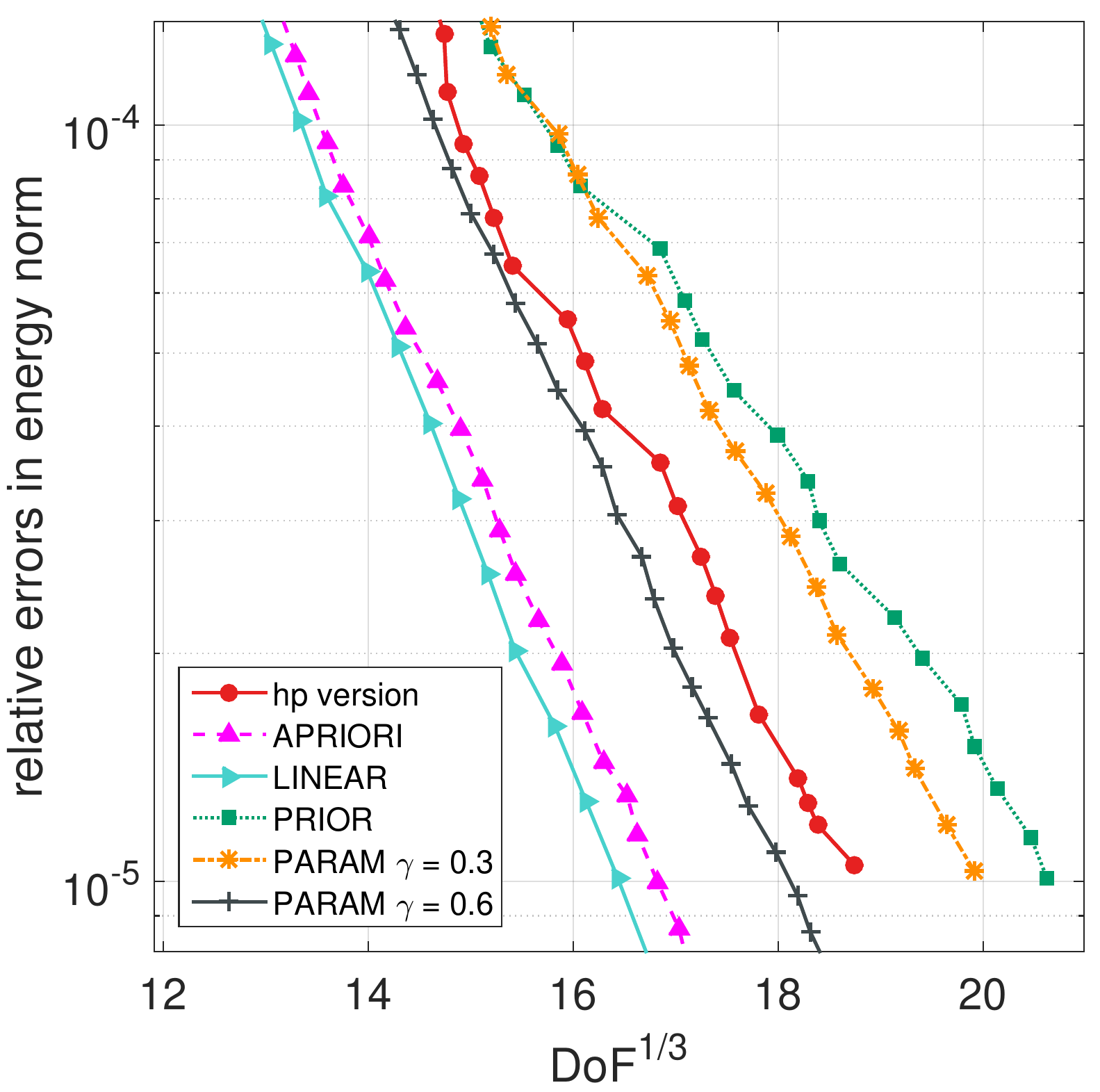}
\end{minipage}
	\caption{
    [L-shape domain of Section~\ref{sec_sing}] Relative energy error $\norm{\Gr (\uu - \uh)}/\norm{\Gr \uu}$
     as a function of $\mathrm{DoF}_\ell^{\frac13}$,
    obtained using the present $hp$-decision criterion, the criteria \textsf{PRIOR} and \textsf{PARAM}
    ($\gamma = 0.3$ and $\gamma = 0.6$), the \textsf{\pd{APRIORI}}, and \textsf{LINEAR} strategy (\textit{left}) and a detailed view (\textit{right}).}
	\label{fig:conv_Lshape}
\end{figure}

\begin{figure}[!ht]
\centering
\begin{minipage}[b]{0.41\linewidth}
	\centering
	\includegraphics[width=0.9\textwidth]{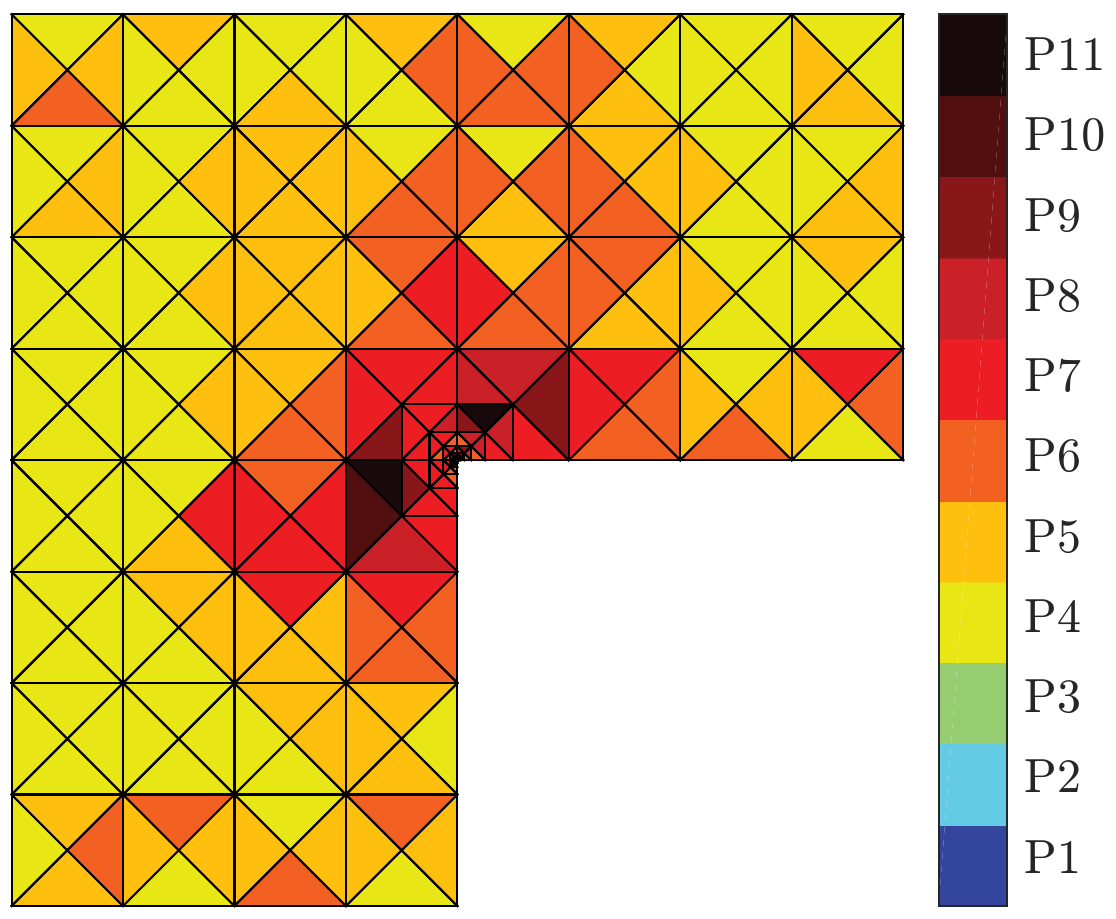}
\end{minipage}
\hfill
\begin{minipage}[b]{0.41\linewidth}
	\centering
	\includegraphics[width=0.76\textwidth]{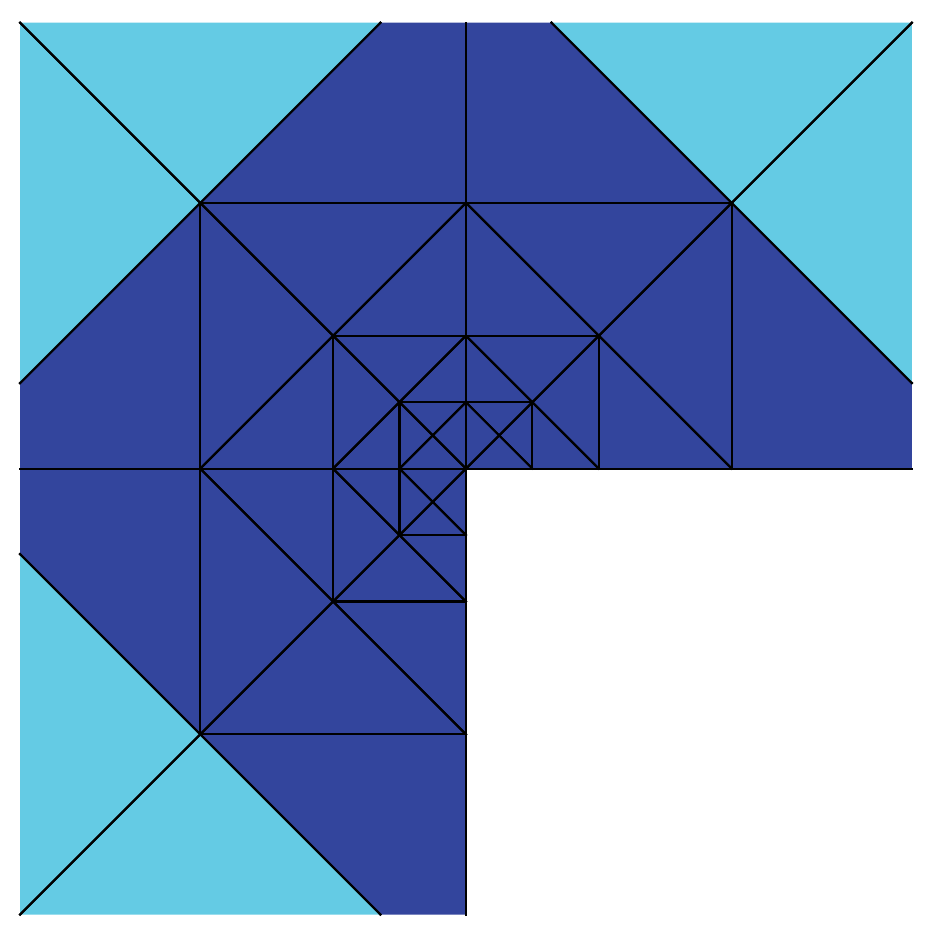}
\end{minipage}
	\caption{[L-shape domain of Section~\ref{sec_sing}] Mesh and polynomial-degree distribution obtained after $70$ iterations (when the relative error reaches $10^{-5}$) of the adaptive procedure employing the \textsf{\pd{APRIORI}} $hp$-strategy (\textit{left}) and a zoom in $[-10^{-7},10^{-7}] \times [-10^{-7},10^{-7}]$ near the re-entrant corner (\textit{right}).
    }
	\label{fig:final_mesh_apriori}
\end{figure}

\begin{figure}[!ht]
\centering
\begin{minipage}[b]{0.41\linewidth}
	\centering
	\includegraphics[width=0.91\textwidth]{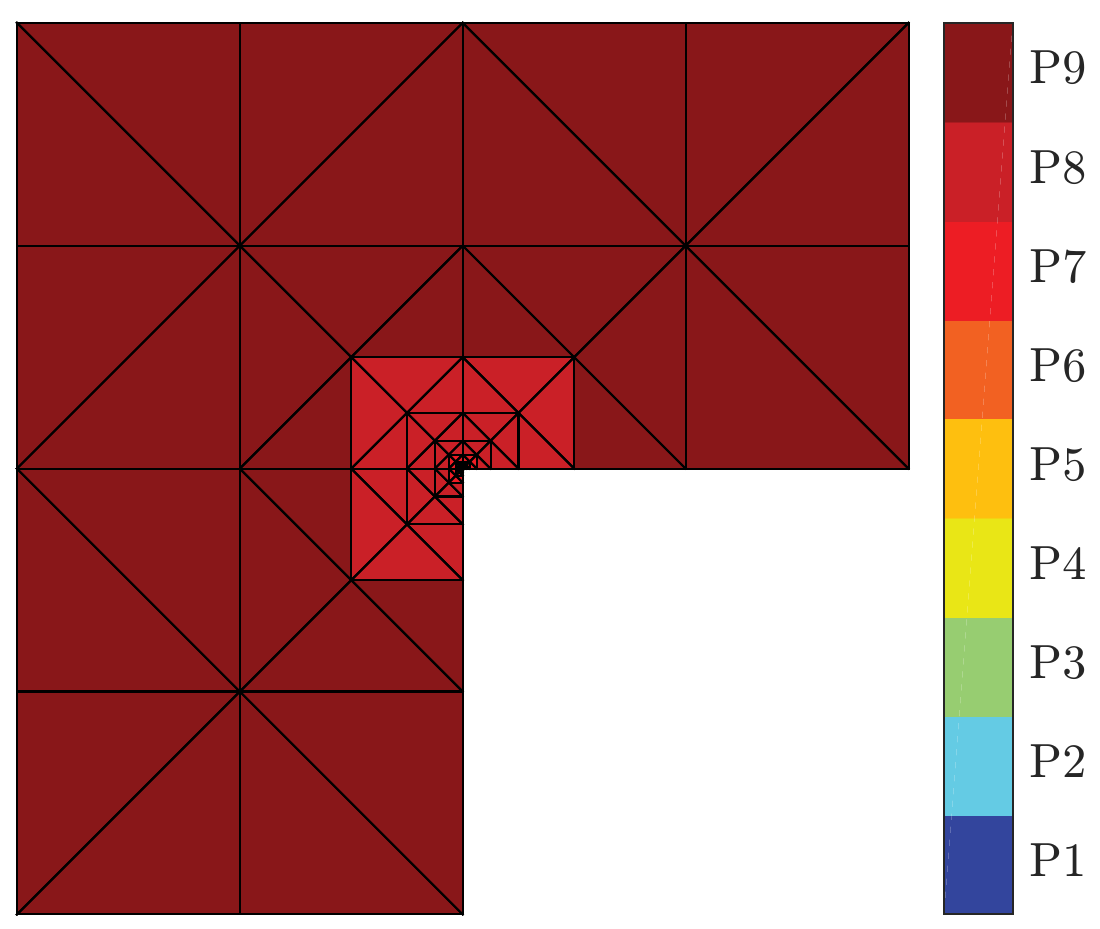}
\end{minipage}
\hfill
\begin{minipage}[b]{0.41\linewidth}
	\centering
	\includegraphics[width=0.75\textwidth]{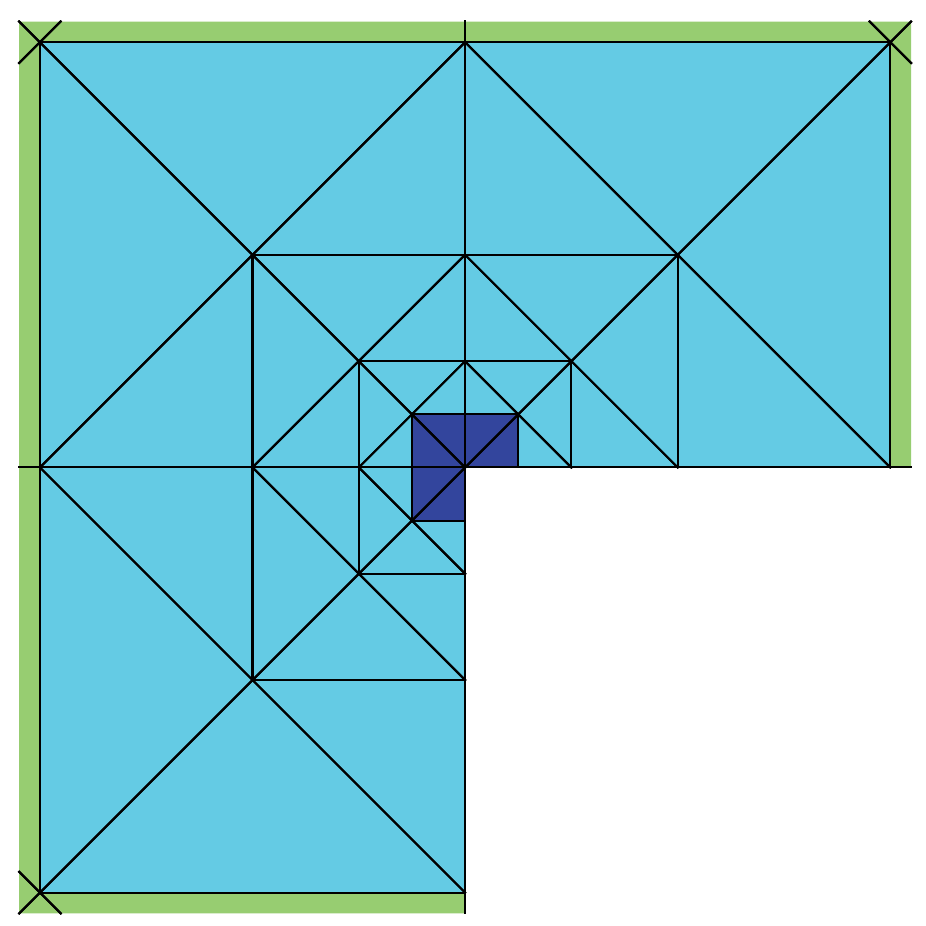}
\end{minipage}
	\caption{[L-shape domain of Section~\ref{sec_sing}] Mesh and polynomial-degree distribution obtained after $45$ iterations (when the relative error reaches $10^{-5}$) of the procedure employing the refinement strategy \textsf{LINEAR} (\textit{left}) and a zoom in $[-10^{-6},10^{-6}] \times [-10^{-6},10^{-6}]$ near the re-entrant corner (\textit{right}).
    }
	\label{fig:final_mesh_LINEAR}
\end{figure}

\begin{figure}[!ht]
	\centering
	\includegraphics[width=0.4\textwidth]{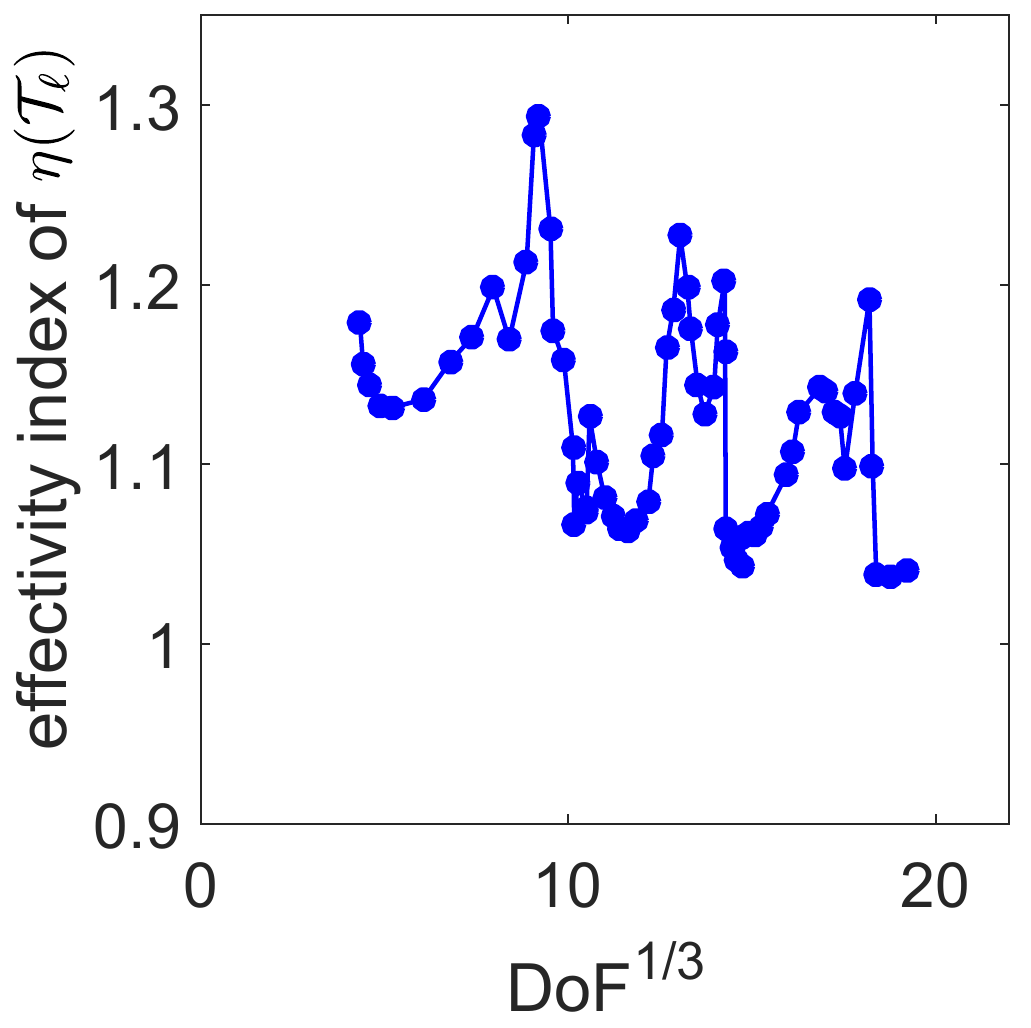}
	\caption{[L-shape domain of Section~\ref{sec_sing}]
	    The effectivity indices of the error estimate $\eta(\Th)$, defined as $\eta(\Th) / \norm{\Gr(\uu - \uh)}$, throughout the 65 iterations of the present $hp$-adaptive procedure.}
	\label{fig:Ieff_estim_lshape}
\end{figure}

\begin{figure}[!ht]
\begin{minipage}[b]{0.49\linewidth}
	\centering
	\includegraphics[width=0.84\textwidth]{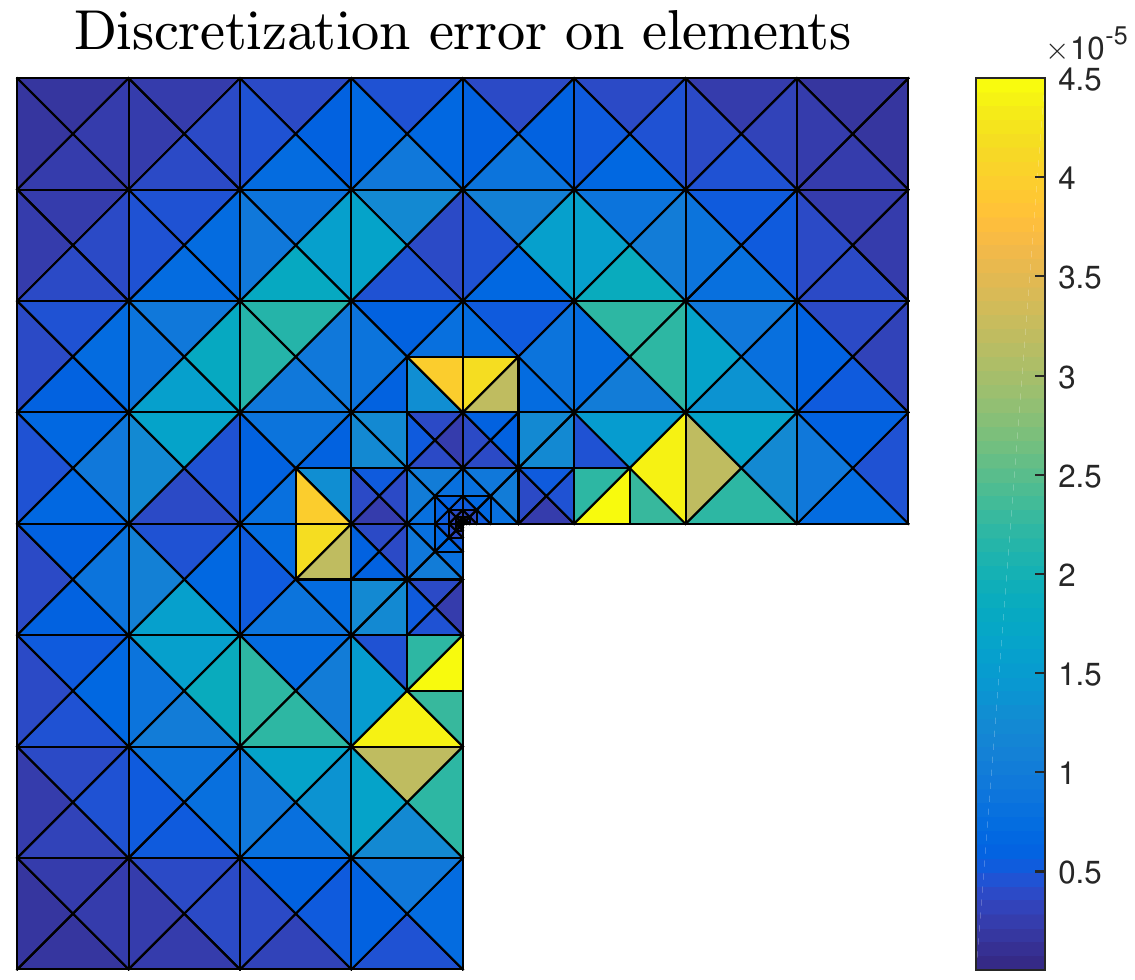}
\end{minipage}
\hfill
\begin{minipage}[b]{0.49\linewidth}
	\centering
	\includegraphics[width=0.87\textwidth]{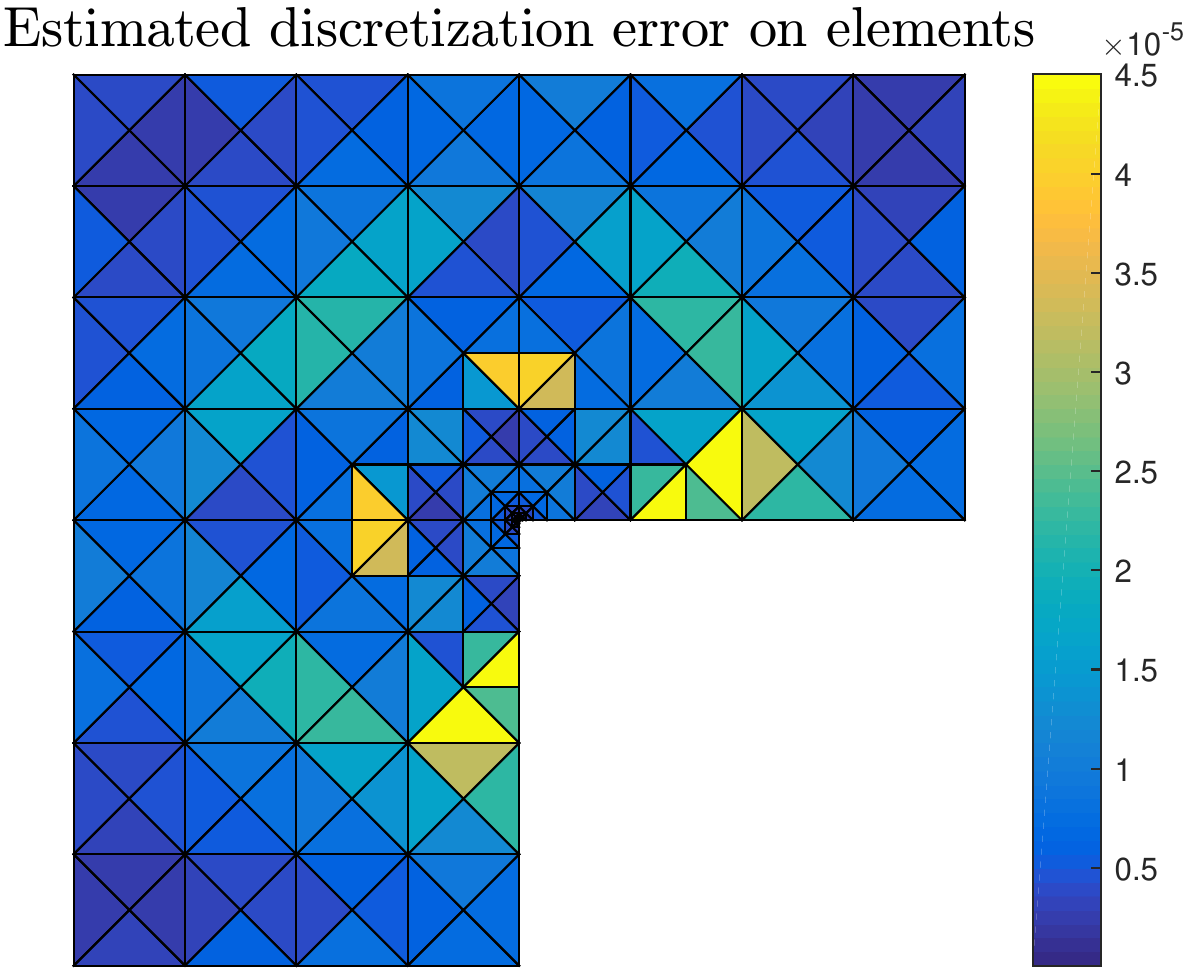}
\end{minipage}
	\caption{[L-shape domain of Section~\ref{sec_sing}] Distribution
		of the energy error $\norm{\Gr(\uu - \uh)}_\elm$ (\textit{left})
	    and of the local error estimators $\eta_\elm$ from Theorem~\ref{thm_est}
	    (\textit{right}), $\ell = 45$.
	    The effectivity index of the estimate defined as $\eta({\T_{45}}) / \norm{\Gr(\uu - \uu_{45})}$ is $1.0468$.}
	\label{fig:estimator_iter_5}
\end{figure}

\begin{figure}[!ht]
\begin{minipage}[b]{0.49\linewidth}
	\centering
	\includegraphics[width=0.82\textwidth]{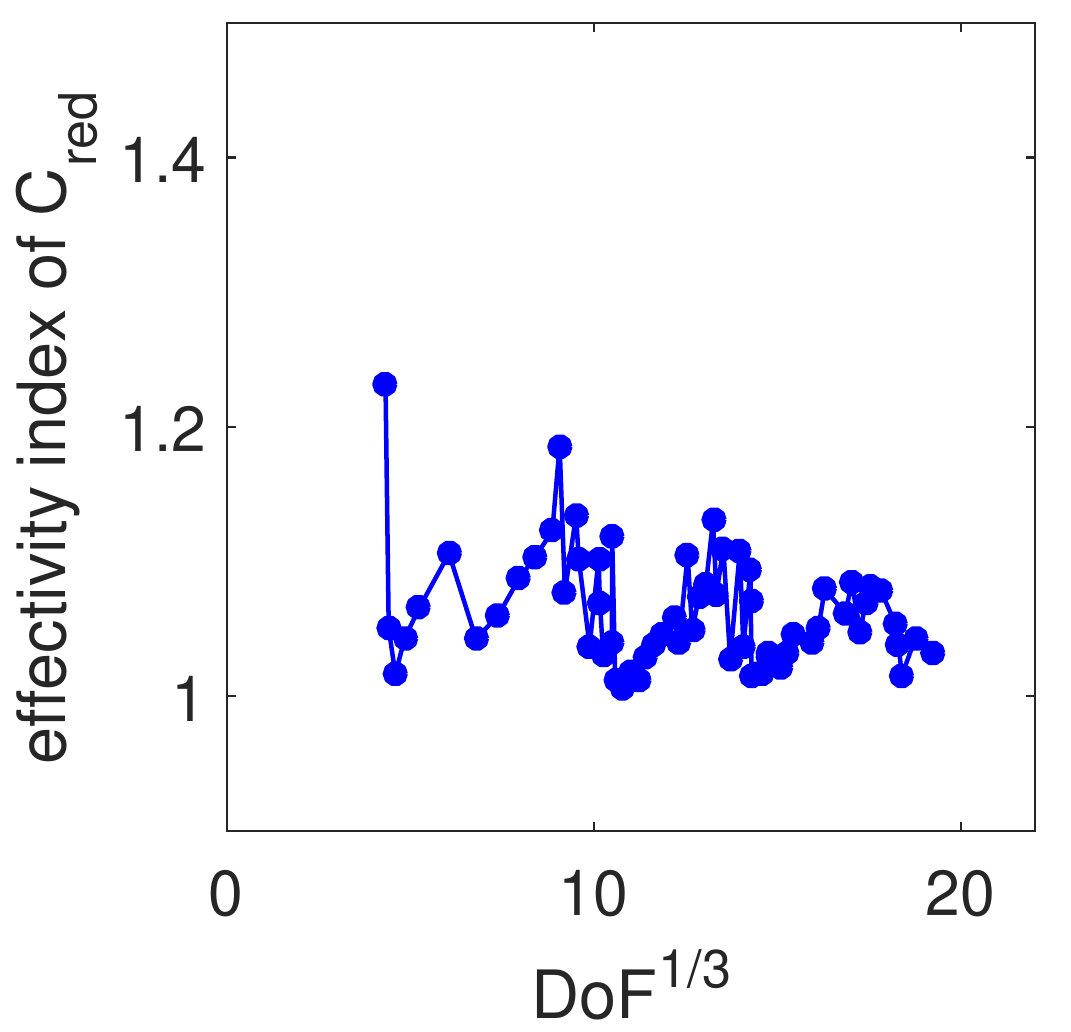}
\end{minipage}
\hfill
\begin{minipage}[b]{0.49\linewidth}
	\centering
	\includegraphics[width=0.80\textwidth]{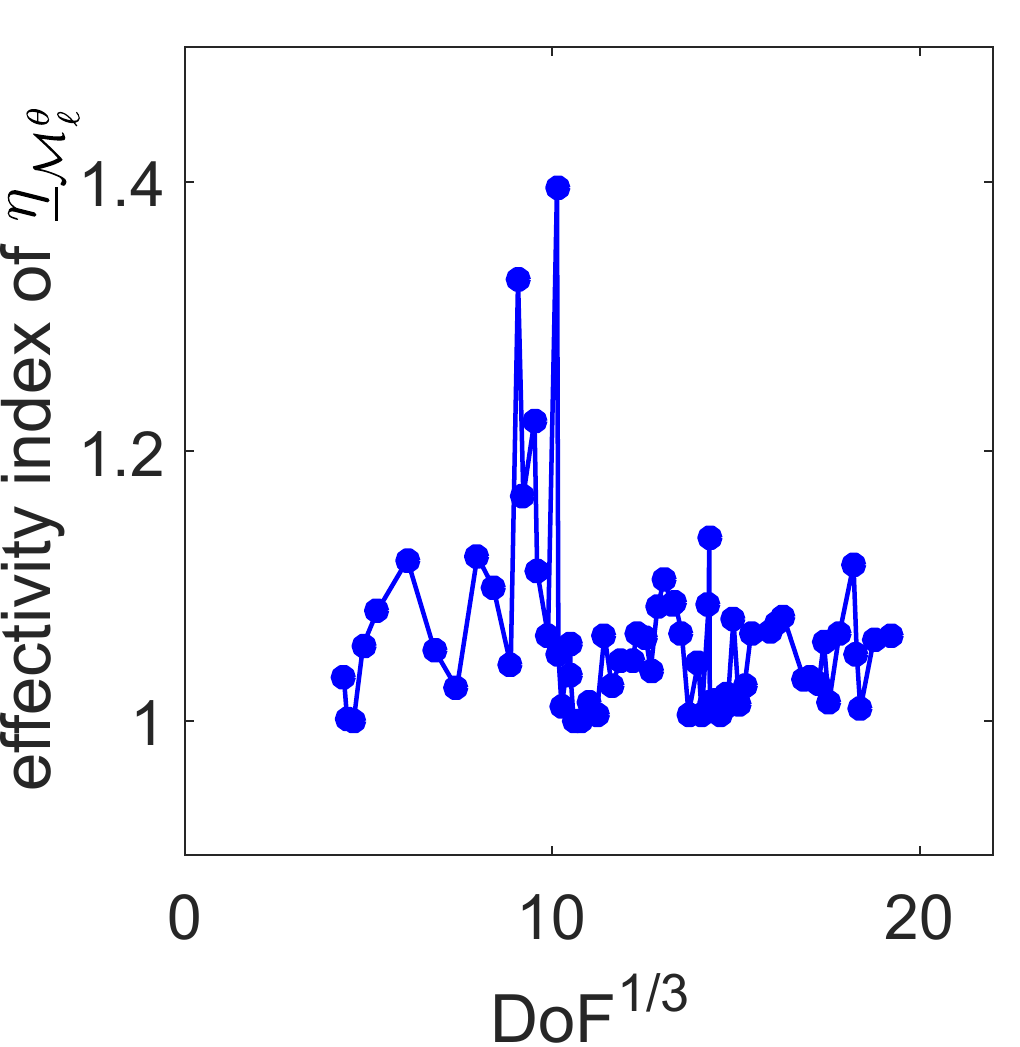}
\end{minipage}
	\caption{[L-shape domain of Section~\ref{sec_sing}] Effectivity indices~\eqref{eq_eff_index} for the error reduction factor $\Cred$ from Theorem~\ref{th_contraction} (\textit{left}) and effectivity indices for the lower bound $\underline \eta_{\Mell}$ from Lemma~\ref{lem_difference_bound_with_res} defined as the ratio $\norm{\Gr(\uhh~-~\uh)}_\omel~/~\underline \eta_{\Mell}$ (\textit{right}).}
	\label{fig:Ccontr_Lshape}
\end{figure}

Figure~\ref{fig:Final mesh} presents the final mesh and polynomial-degree
distribution after $65$ steps of the $hp$-adaptive
procedure~\eqref{eq_afem_loop} (left panel) along with a zoom in the window
$[-10^{-6},10^{-6}] \times [-10^{-6},10^{-6}]$ near the re-entrant corner
(right panel). Figure~\ref{fig:conv_Lshape} (left panel) displays the
relative error $\norm{\Gr(\uu-\uh)}/\norm{\Gr \uu}$ as a function of
$\mathrm{DoF}_\ell^{1/3}$ in logarithmic-linear scale to illustrate that, as
in the previous test case, the present $hp$-adaptive procedure leads to an asymptotic exponential rate of convergence.
The corresponding values of constants~$C_{1}$ and~$C_{2}$ in expression~\eqref{eq_exp_converg} obtained by the 2-parameter least squares fit are \num[round-mode=places,round-precision=2]{4.73153} and \num[round-mode=places,round-precision=2]{0.686994}, respectively.
\pd{For the direct comparison with other methods, we refer to the long version~\cite[Table~15]{Mitchell_hp_review_long} of the survey paper~\cite{Mitchell_hp_review_14}.}
However, note that data sets of greater sizes than in our case have been used for the least squares fitting therein.
A detailed view when the error
takes lower values is provided in the right panel of
Figure~\ref{fig:conv_Lshape}. We also plot the relative
errors obtained when using the $hp$-decision criteria \textsf{PRIOR} and
\textsf{PARAM}, as well as those obtained using \pd{\textsf{\pd{APRIORI}}} criterion 
exploiting the a priori knowledge of the exact solution (marked simplices are $h$-refined
only if they touch the corner singularity, otherwise they are $p$-refined).
\pd{In addition, we provide also the relative errors obtained by employing the 
   (non-adaptive) strategy which we refer to as \textsf{LINEAR}, inspired by the theoretical results for the one-dimensional problem with singular solution~\cite{Gui_Bab_hp_II_86,Gui_Bab_hp_III_86,Szab_Bab_FEA_91}.
   When employing this strategy, we start from a coarse grid $\T_0$ with $\max_{\elm \in
   \T_0}h_\elm = 0.5$.
   At each iteration, only the patch containing the re-entrant corner is $h$-refined. 
   Thus, the elements of each mesh $\Th$, $\ell \ge 1$, decrease in size in geometric progression (in our case with factor 0.5) toward the re-entrant corner.
   For each $\Th$, $\ell \ge 1$, we group the elements in layers $\mathcal{L}_{1},\mathcal{L}_{2},\ldots, \mathcal{L}_{m(\ell)}$  depending on their distance from the origin ($\mathcal{L}_{1}$ containing the singularity), such that $\Th = \bigcup_{i = 1}^{m(\ell)} \mathcal{L}_{i}$. The total number of layers $m(\ell)$ depends on how many times the current mesh $\Th$ has been refined. 
   Each element $\elm \in \Th$ is then assigned polynomial degree $p_{\ell,\elm}$ layer-wise, increasing linearly away from the singularity, in the way
   $$ p_{\ell,\elm} \eq \left\lceil 1 + \frac{(i-1)}{3} \right\rceil,$$
   where $i$ is the index of the layer $\mathcal{L}_{i}$ containing the element $\elm$.}   
\pd{For the strategy \textsf{\pd{APRIORI}} (Figure~\ref{fig:final_mesh_apriori})  and the strategy \textsf{LINEAR} (Figure~\ref{fig:final_mesh_LINEAR}), we illustrate 
also the resulting polynomial-degree distribution at 
the step when the relative error reaches $10^{-5}$.
As for the previous test case, in Figure~\ref{fig:Ieff_estim_lshape} we illustrate the quality of the error estimator from Theorem~\ref{thm_est} in terms of the effectivity
 index $\eta(\Th)/\norm{\Gr (\uu - \uh)}$ throughout all the iterations of the present $hp$-adaptive process. 
 Figure~\ref{fig:estimator_iter_5} then compares the actual and estimated error
 distributions on iteration $\ell=45$ of the adaptive loop, showing excellent
 agreement.}
Figure~\ref{fig:Ccontr_Lshape} (left panel) presents the
effectivity index for the reduction factor $\Cred$, see~\eqref{eq_eff_index},
throughout the adaptive process, whereas the right panel of
Figure~\ref{fig:Ccontr_Lshape} examines the quality of the lower bound
$\underline \eta_{\Mell}$ from Lemma~\ref{lem_difference_bound_with_res} by
plotting the ratio of the left-hand side to the right-hand side of the lower bound
in~\eqref{eq_res_norm_low_bound}. For both quantities, we can
draw similar conclusions to the previous test case, thereby
confirming that sharp estimates on the error reduction factor are available.
Additional numerical experiments (not shown here) indicate that the lower
bound estimate can be made even sharper by performing $h$-refinement so as to
satisfy the interior node property. Finally, to give some further insight
into the $hp$-adaptive process, we present in Tables~\ref{table:history_1}
and~\ref{table:history_2} some details on the $hp$-refinement decisions made
by the proposed $hp$-refinement criterion during the first 10 and the last 10
iterations of the adaptive loop. We observe that in the initial iterations,
where the underlying mesh is still rather coarse, the polynomial degree is
increased also on those simplices touching  the re-entrant corner.
Nevertheless, this decision does not occur anymore later when the mesh around
the singularity is already more strongly refined than in the rest of the
domain. Therefore, an improvement of our approach is expected, as suggested
in \cite{Can_Noch_Stev_Ver_hp_AFEM_16}, in conjunction with an appropriate
coarsening strategy correcting the excessive $p$-refinement in the early
stages.
\pd{Table~\ref{table:comp_iter_dofs} (bottom) again brings some additional 
comparisons with other strategies in terms of number of iterations and
number of degrees of freedom necessary to reach relative error $10^{-5}$.
We observe that the results achieved using the 
present strategy are comparable with those achieved by other (established) strategies.}

\begin{table}[h!]
\centering
\begin{tabular}{l cccccccccc }
 \toprule
 Iteration & 1 & 2 & 3 & 4 & 5 & 6 & 7 & 8 & 9 & 10 \\ \hline
 Triangles  & 192 & 192 & 192 & 192 & 192 & 198 & 204 & 210 & 216 & 222\\
 Maximal polynomial degree  & 1 & 2 & 3 & 4 & 5 & 5 & 5 & 5 & 5 & 5 \\
 Marked vertices   & 1 & 1 & 1 & 2 & 2 & 3 & 4 & 4 & 6 & 6\\
 Triangles flagged for $h$-refinement  & 0 & 0 & 0 & 0 & 6 & 6 & 6 & 6 & 6 & 6 \\
 Triangles flagged for $p$-refinement & 6 & 6 & 6 & 12 & 6 & 12 & 16 & 18 & 16 & 18 \\
 Triangles flagged for $hp$-refinement &0 & 0 & 0 & 0 & 2 & 0 & 0 & 0 & 0 & 0  \\
 \hline
 \bottomrule
\end{tabular}
\caption{[L-shape domain of Section~\ref{sec_sing}] Refinement decisions in Algorithm~\ref{hp_strat_alg} during the first 10 iterations of the adaptive loop~\eqref{eq_afem_loop}.
}
\label{table:history_1}
\end{table}

\begin{table}[h!]
\centering
\begin{tabular}{l cccccccccc }
 \toprule
 Iteration &  56 & 57 & 58 & 59 & 60 & 61 & 62 & 63 & 64 & 65 \\ \hline
 Triangles  & 492 & 512 & 518 & 524 & 538 & 568 & 574 & 580 & 614 & 660 \\
 Maximal polynomial degree  &  5	& 5	& 5	& 5	& 5	& 5	& 5	& 5	& 5	& 5 \\
 Marked vertices   & 4 & 4 & 5 & 4 & 3 & 3 & 3 & 4 & 5 & 5 \\
 Triangles flagged for $h$-refinement  &  16 & 6 & 6 & 6 & 18 & 6 & 6 & 28 & 30 & 32  \\
 Triangles flagged for $p$-refinement &  8 & 16 & 13 & 22 & 0 & 6 & 6 & 0 & 0 & 8  \\
 Triangles flagged for $hp$-refinement &  0 & 0 & 0 & 0 & 0 & 0 & 0 & 0 & 0 & 4  \\
 \hline
 \bottomrule
\end{tabular}
\caption{[L-shape domain of Section~\ref{sec_sing}] Refinement decisions in Algorithm~\ref{hp_strat_alg} during the last 10 iterations of the adaptive loop~\eqref{eq_afem_loop}.}
\label{table:history_2}
\end{table}

\begin{table}[]
\centering
\begin{tabular}{l|l|c|c|c|c|c}
	           &             		& our  & \pd{APRIORI} & PRIOR   & PARAM 0.3 & PARAM 0.6 \\ \hline
Sharp Gaussian & iter        		& 27           & --       & 37      & 36        & 40        \\
(relative error $10^{-3}$)               & $\text{DoF}^{1/3}$ & 12.56      & --       & 14.29 & 14.06   & 12.49   \\ \hline
L-shape domain & iter        		& 65           & 70       & 68      & 67        & 68        \\
(relative error $10^{-5}$)               & $\text{DoF}^{1/3}$ & 19.24     & 17.35  & 20.82 & 20.07   & 18.18
\end{tabular}
\caption{Comparison of the different adaptive $hp$-strategies in terms of the number of iterations of the loop~\eqref{eq_afem_loop} and of the number of degrees of freedom necessary to reach the given relative error for model problems of Sections~\ref{sec_reg} and~\ref{sec_sing}.}
\label{table:comp_iter_dofs}
\end{table}
\FloatBarrier

\section{Conclusions}
\label{sec:conc}

In this work, we have devised an $hp$-adaptive strategy to approximate model
elliptic problems using conforming finite elements.
Mesh vertices are marked using polynomial-degree-robust a
posteriori error estimates based on equilibrated fluxes. Then marked vertices
are flagged either for $h$- or for $p$-refinement based on the solution of
two local finite element problems where local residual liftings are computed.
Moreover, by solving a third local finite element problem once the
$hp$-decision has been taken and the next mesh and polynomial-degree
distribution have been determined, it is possible to compute a
guaranteed bound on the error reduction factor. Our numerical experiments featuring
two-dimensional smooth and singular weak solutions indicate that the present
$hp$-adaptive strategy leads to asymptotic exponential convergence rates
with respect to the total number of degrees of freedom employed to compute
the discrete solution. Moreover, our bound on the error reduction factor
appears to be, in most cases, quite sharp. Several extensions of the present
work can be considered. On the theoretical side, it is important to prove
that our bound on the reduction factor $\Cred$ is smaller than one and to
study how it depends on the mesh-size and especially on the polynomial
degree. On the numerical side, three-dimensional test cases and taking into
account an inexact algebraic solver are on the agenda.
\FloatBarrier


\def\polhk#1{\setbox0=\hbox{#1}{\ooalign{\hidewidth
  \lower1.5ex\hbox{`}\hidewidth\crcr\unhbox0}}} \def\cprime{$'$}

\end{document}